    %%%%%%%%%%%%%%%%%%%%%%%%%%%%%%%%%%%%%%%%%%%%%%%%%%%%  
    %%                                                %%  
    %%  BRAIDING STRUCTURES ON FORMAL POISSON GROUPS  %%  
    %%        AND CLASSICAL SOLUTIONS OF THE          %%  
    %%         QUANTUM YANG-BAXTER EQUATION           %%  
    %%                                                %%  
    %%%%%%%%%%%%%%%%%%%%%%%%%%%%%%%%%%%%%%%%%%%%%%%%%%%%  

\input amstex
\documentstyle{amsppt}

\magnification=\magstep1  
\hsize=6.5truein  
\vsize=9truein 

%\baselineskip=.15truein
\baselineskip=12pt 

\loadbold 

  % at 10truept 
  % at 7truept 
  % at 10truept 
  % at 10truept 
  % at 10truept 
% \font \eusb=Eusb10.mf   % at 10truept 
\font \eusb=cmtt10  at 11pt 

\def \loongrightarrow {\relbar\joinrel\relbar\joinrel\rightarrow}
\def \llongrightarrow 
   {\relbar\joinrel\relbar\joinrel\relbar\joinrel\rightarrow} 
 
\def \gerg {\frak g}

\def \gersl {\frak{sl}}
\def \e {\hbox{\rom{e}}}
\def \f {\hbox{\rom{f}}}
\def \h {\hbox{\rom{h}}}
\def \r {\boldkey{r}}
\def \N {\Bbb N}

\def \O {\Cal O} 
\def \R {\frak R} 
\def \lagR {\hbox{\eusb R}} 
\def \HA {{\Cal H}{\Cal A}}  
\def \Ker {\hbox{\rm Ker}}
\def \kh {\Bbbk[[\hbar]]}
\def \kqqm {\Bbbk\big[q,q^{-1}\big]}
\def \ug {U(\gerg)}
\def \uhg {U_\hbar(\gerg)}
\def \uhat {\widehat{U}}

\def \utilde {\widetilde{U}}

\def \fgstar {F[[\gerg^*]]}

\def \otimeshat {\,\widehat{\otimes}\,}
\def \otimestilde {\,\widetilde{\otimes}\,}

\def \QUEA {\hbox{$ \displaystyle{\Cal{Q\hskip1ptUE\hskip-1ptA}} $}}
\def \QTQUEA {\hbox{  %  
   $ \displaystyle{\Cal{QT\hskip-1ptQ\hskip1ptUE\hskip-1ptA}} $}} 
\def \TQUEA {\hbox{   %  
   $ \displaystyle{\Cal{T\hskip-1ptQ\hskip1ptUE\hskip-1ptA}} $}} 
\def \QFSHA {\hbox{$ \displaystyle{\Cal{Q\hskip1ptFSHA}} $}}
\def \BQFSHA {\hbox{$ \displaystyle{\Cal{BQF\hskip1ptSHA}} $}} 
\def \RBQFSHA {\hbox{$ \displaystyle{\Cal{RBQ\hskip1ptFSHA}} $}} 
\def \lsF {\hbox{$ {\phantom{\big|}}^{\scriptscriptstyle F\!\!} $}}

\document

\topmatter 

\title  
  Braiding structures on formal Poisson groups  \\   
  and classical solutions of the QYBE  
\endtitle

\author
    Fabio Gavarini${}^\dag$  \;  and  \;  Gilles Halbout${}^\ddag$   
\endauthor 

% \leftheadtext{ Fabio Gavarini, \ Gilles Halbout }   
% \leftheadtext{ Fabio Gavarini  {\ }  and  {\ }  Gilles Halbout }   
\leftheadtext{ Fabio Gavarini  \;  and  \;  Gilles Halbout }   

\rightheadtext{ Braiding structures and classical solutions of the 
QYBE }   

\affil 
  ${}^\dag$ \! Universit\`a di Roma ``Tor Vergata''  \\   
Dipartimento di Matematica  --  Roma, ITALY  \\ 
  {}  \\  
  ${}^\ddag$ \! \hbox{Institut de Recherche
Math\'ematique Avanc\'ee,}  \\   
  ULP--CNRS   --   Strasbourg, FRANCE  \\ 
\endaffil

\address\hskip-\parindent
  ${}^\dag$ \! Universit\`a degli Studi di Roma ``Tor Vergata''  
---   Dipartimento di Matematica   \newline
  Via della Ricerca Scientifica, 1   ---   I-00133 Roma,
ITALY  \newline
  e-mail: \ gavarini\@{}mat.uniroma2.it  \newline
     \newline  
  ${}^\ddag$ \! Institut de Recherche Math\'ematique Avanc\'ee  
---   Universit\'e{} "Louis Pasteur" -- C.N.R.S.  \newline 
  7, rue Ren\'e{} Descartes   ---   F-67084 STRASBOURG Cedex,
FRANCE  \newline
  e-mail: \  halbout\@math.u-strasbg.fr  
\endaddress

\abstract  
   If  $ \gerg $  is a quasitriangular Lie bialgebra, the formal
Poisson group  $ \fgstar $  can be given a braiding structure: 
this was achieved by Weinstein and Xu using purely geometrical 
means, and independently by the authors by means of quantum 
groups. In this paper we compare these two approaches: first, we 
show that the braidings they produce share several similar 
properties (in particular, the construction is functorial); 
second, in the simplest case ($ \, G = {SL}_2 \, $)  they do 
coincide.  The question then rises of whether they are always the 
same.   
\endabstract  

% 
% \keywords  
%    Quasitriangular Poisson groups, quantum groups, 
% quantum Yang-Baxter equation  
% \endkeywords  
% 
% 
% \subjclass  
%    Primary 17B37, 20G42, Secondary 81R50, 16W30  
% \endsubjclass  
% 

\endtopmatter

   \footnote""{Keywords: \ {\sl Quantum Groups, quasitriangular
Poisson groups, quantum Yang-Baxter equation}.}  
   \footnote""{2000  {\it Mathematics Subject Classification:} 
\,  Primary 17B37, 20G42, Secondary 81R50, 16W30. }   
%    \footnote""{Revised version of September 
% $ 4^{\text{th}}$,  2002.}  

\vskip20pt

\centerline {\bf Introduction } 

\vskip13pt 

   In the study of classical Hamiltonian systems, one is naturally
interested in those which are completely integrable.  A natural
condition to achieve complete integrability for the system is that
it admit a so called "Lax pair", thus one typical goal is to find
Hamiltonian systems admitting such a pair; a standard recipe to
obtain this has been provided by Semenov-Tian-Shansky (see [Se]),
which explain how to get such a system proceeding from a pair 
$ (\gerg,\r) $  where  $ \gerg $  is a Lie quasitriangular Lie
bialgebra and  $ \r $  is its  $ \r $--matrix,  a classical
solution of the classical Yang-Baxter equation (CYBE): the
system is built up on  $ \gerg^* $,  the Lie bialgebra dual 
to  $ \gerg $,  as phase space, and the  $ \r $--matrix  $ \r $ 
provides (a recipe for) the Poisson bracket on  $ C^\infty
(\gerg^*) $.  This raises the question of studying quasitriangular
bialgebras, as objects of special interest within the category of
Lie bialgebras: in particular, since we think at  $ \gerg^* $  as
a phase space, so that  $ \gerg $  is its cotangent space, one's
desire is to understand the geometrical meaning of the classical 
$ \r $--matrix.   
                                             \par  
   A second motivation for studying the geometrical meaning of
the classical  $ \r $--matrix  arises from conformal, quantum
and topological quantum field theories.  Indeed, all these are
concerned with the notion of "fusion rules" which, roughly,
rule the tensor product in a quasitensor category (see e.g.~[FK]):
as an application   --- among others ---   one has a recipe which
provides tangle and link invariants as well as invariants of
3-manifolds (cf.~[Tu]).  In this setting, the common notion
one start with is that of a quasitensor (or "braided monoidal")
category; such an object can be built up as category of
representations of a quasitriangular Hopf algebra (QTHA):
indeed, by Tannaka-Krein reconstruction theorems the two notions  
--- quasitensor categories and quasitriangular Hopf algebras
---   are essentially equivalent, so one may switch to the study
of QTHAs.  A key example of QTHA is given by a quantum group, in
the shape of a quantum universal enveloping algebra (QUEA) together
with its (universal)  $ R $--matrix.  Now, the semiclassical
counterpart of a QUEA is a Lie bialgebra  $ \gerg $  (i.e., the
given QUEA is the quantization of  $ U(\gerg) $):  if the QUEA is
also quasitriangular, then the semiclassical counterpart of its 
$ R $--matrix  is a classical  $ \r $--matrix  $ \r $  on  $ \gerg $, 
the pair  $ (\gerg,\r) $  being a quasitriangular Lie bialgebra. 
The question then rises of whether   --- or at least how far ---  
one can perform the constructions which are usually made via the
QUEA and its  $ R $--matrix  (such as that of link invariants)
using instead only the "semiclassical" datum of  $ (\gerg,\r) $: 
then again the key point will be to understand the geometrical
meaning of the classical  $ \r $--matrix.   
                                             \par  
   With this kind of motivations, we go and study the following
problem.  It is known that if  $ \gerg $  is a Lie bialgebra (over
a field  $ \Bbbk $  of zero characteristic), then its dual space 
$ \gerg^* $  is a Lie bialgebra as well.  Also, let  $ G $  be an
algebraic Poisson group   --- or Poisson-Lie group, say, when  $ \,
\Bbbk \in \{\Bbb{R,C}\} \, $ ---   whose tangent Lie bialgebra is 
$ \gerg $.  Now assume  $ \gerg $  is quasitriangular, with 
$ \r $--matrix  $ \r \, $:  this gives to $ \gerg $  some
additional properties; two questions then rise:   
                                             \par  
   ($ * $) \; What an additional structure  one obtains on the
dual Lie bialgebra  $ \gerg^* \, $?  
                                             \par  
   ($ \bullet $) \; What is the geometrical global datum on 
$ G $  which is the result of "integrating"  $ \r \, $?  
                                             \par  
   Of course, the two questions and their answers are necessarily
tightly related.  
                                             \par  
   First, an answer to question  $ (*) $  was given by the authors
in [GH] (cf.~also [Re], [Ga1], [Ga2]): the topological Poisson Hopf 
algebra  $ \fgstar $  (the function algebra of the formal Poisson
group associated to  $ \gerg^* $) is  {\sl braided\/}  (see the
definition later on).  
                                             \par  
  The result in [GH] was proved using the theory of quantum groups. 
Indeed, after Etingof-Kazhdan (cf.~[EK]) every Lie bialgebra admits
a quantization  $ \uhg $,  namely a (topological) Hopf algebra over 
$ \kh $  whose specialisation at  $ \, \hbar = 0 \, $  is isomorphic
to  $ \, \ug \, $  as a co-Poisson Hopf algebra; in addition, if 
$ \gerg $  is quasitriangular and  $ \r $  is its  $ \r $--matrix, 
then such a  $ \uhg $  exists which is quasitriangular too, as a
Hopf algebra, with an  $ R $--matrix  $ \, R_\hbar \, (\, \in \uhg
\otimes \uhg \, $)  such that  $ \; R_\hbar \equiv 1 + \r \, \hbar
\; \mod\, \hbar^2 \; $  (here one identifies, as  $ \kh $--modules, 
$ \, \uhg \cong \ug [[\hbar]] \, $).  Using Drinfeld's  {\it Quantum
Duality Principle}  ([Dr1]; cf.~[Ga5] for a proof), from any QUEA 
$ \uhg $  with semiclassical limit  $ \ug $  one can extract a
certain quantum formal series Hopf algebra (QFSHA)  $ {\uhg}' $ 
such that the semiclassical limit of  $ {\uhg}' $  is  $ \fgstar $. 
In [GH], starting from a quasitriangular QUEA  $ \big( \uhg, R \big) $, 
we showed that, although  {\sl a priori\/}  $ \, R \not\in {\uhg}'
\otimes {\uhg}'  
%% 
%%  = {\big( \uhg \otimes \uhg \big)}'  
%%  
\, $  (so that the pair  $ \, \big( {\uhg}', R \,\big) \, $  is 
{\sl not}  in general a quasitriangular Hopf algebra), nevertheless 
its adjoint action  
   \hbox{$ \; \R_\hbar := {\hbox{\rm Ad}}(R_\hbar) : \,
\uhg \!\otimes\! \uhg \longrightarrow \uhg \!\otimes\! \uhg \, , 
\;  x \!\otimes\! y \mapsto R_\hbar \!\cdot\! (x \!\otimes\! y)
\!\cdot\! R_\hbar^{\,-1} $}   
\; stabilises the subalgebra  $ \, {\uhg}' \otimes {\uhg}' \, $, 
\, hence induces by specialisation an operator  $ \R_0 $ over  $ \, 
\fgstar \otimes \fgstar \, $:  moreover, the properties which make 
$ R_\hbar $  an  $ R $--matrix  imply that  $ \R_\hbar $ is a 
braiding operator, hence the same holds for  $ \R_0 $:  thus, the 
pair  $ \big(\fgstar, \R_0\big) $  is a braided Hopf algebra.  In 
particular,  this gives us a new method to produce set-theoretical 
solutions of  the QYBE, thus giving a positive answer to a 
question set in [Dr2]  (also tackled, for instance, in [ESS]).  
%%  
%%  ; in 
%%  addition, this yields at once representations of the braid groups, 
%%  which is the starting point for the construction of knot and link 
%%  invariants following the pattern in [Tu] (cf.~also [CP]).  
%%  
   Note also that for igniting  our construction we only need
a quantisation functor  $ \; (\gerg,\r) \mapsto \big( \uhg, R 
\,\big) \, $,  
%%  
%%  \, which by no means must be the one of [EK];
%%  any other would work, for instance Enriquez' one (see [En]).  
%%  
    and several of them exist (see [En]).   
                                             \par  
   Second, an answer to question  $ (\bullet) $  was given by
Weinstein and Xu in [WX].  We briefly sketch their results.  Let 
$ G $,  resp.~$ G^* $,  be a Poisson group with tangent Lie
bialgebra  $ \gerg $, resp.~$ \gerg^* \, $:  \, in addition,
assume both  $ G $  and  $ G^* $  to be complete.  Let  $ D $ 
be the corresponding double Poisson group, which is given a
structure of symplectic double groupoid, over  $ G $  and 
$ G^* $  at once (further assumptions are needed, see \S 3
later on).  Then the authors prove that there is a classical
analogous of the quantum  $ R $--matrix,  namely a Lagrangian
submanifold  $ \lagR $  of  $ D \times D $,  called  {\it the
(global) classical  $ \lagR $--matrix},  which enjoys much the
same properties of a quantum  $ R $--matrix!  Furthermore, for
any symplectic leaf  $ S $ in  $ G^* $,  this  $ \lagR $  induces
a symplectic automorphism of  $ S \times S $  which in turn at the
level of function algebras yields a braiding for  $ F[S] \, $;  \,
then, as  $ G^* $  is the union of its symplectic leaves, we get
also a braiding on  $ F[G^*] $  and so, via completion, a braiding
on  $ \fgstar $  too.   
                                             \par  
   As a first goal in this paper, we investigate more in depth the
properties of the construction in [GH].  In particular, we show that
the step  $ \; \big( \uhg, R \,\big) \mapsto \big( {\uhg}', \R_\hbar
\,\big) \, $  is functorial and preserves quantisation equivalence. 
Since the initial quantisation step $  \, (\gerg,\r) \mapsto \big(
\uhg, R_\hbar \big) \, $  (provided by [EK], but any other would
work) is functorial, and of course the final specialisation step 
$ \, \big({\uhg}',\R_\hbar\big) \mapsto \big( \fgstar, \R_0\big)
\, $  is trivially functorial, we conclude that the whole
construction  $ \, (\gerg,\r) \mapsto \big(\fgstar, \R_0\big) \, $ 
is functorial too.  Moreover, whenever one has a braiding on 
$ \fgstar $  a so-called  {\it infinitesimal braiding\/} 
$ \overline{\R} $  is defined on the cotangent Lie bialgebra
of  $ \fgstar^{\otimes 2} $,  which is just $ \gerg^{\oplus 2} $: 
if the braiding is the afore mentioned  $ \R_0  $,  we prove that
the infinitesimal braiding  $ \overline{\R}_0 $  is trivial.  
                                             \par  
   As a second goal of the paper, we compare our results with those
of [WX].  First of all, a general fact is worth stressing: the purpose
in [WX]  is to find a geometrical counterpart of the classical 
$ \r $--matrix,  in particular an object which is of global rather
than local nature: to this end, one is forced to impose some
additional requirements from scratch, mainly the existence of
complete Poisson groups  $ G $  and  $ G^* $  with tangent Lie
bialgebras respectively  $ \gerg $  and  $ \gerg^* \, $).  In
contrast, the approach of [GH] sticks to the infinitesimal level:
everything is formulated in terms of Lie bialgebras or formal
Poisson groups.  Therefore, the final output of [WX] is stronger
but requires stronger hypotheses as well.  Nevertheless, the
additional requirements in [WX] are not necessary if we stick to
the infinitesimal setting: indeed, a good deal of the analysis
therein can be carried out as well in local terms   --- just on
germs of Poisson groups ---   so that eventually one ends up with
results which are perfectly comparable with those of [GH].  Thus
we compare the braiding $ \R_{{}_{W\!X}} $  of [WX] with the one
of [GH], call it  $ \R_{{}_{G\!H}} $.  Indeed, one has a theoretical
reason to find strong similarities: namely, the construction in [WX]
is a  {\sl geometric\/}  quantisation of  $ (\gerg,\r) $,  whereas
the one of [GH] passes through  {\sl deformation\/}  quantisation. 
As a matter of fact, first we show that the infinitesimal
braiding  $ \overline{\R}_{{}_{W\!X}} $  is trivial, just like 
$ \overline{\R}_{{}_{G\!H}} $.  Second, when  $ \, \gerg =
\gersl_2 \, $  with the standard  $ \r $--matrix  we prove
via explicit computation that  $ \, \overline{\R}_{{}_{W\!X}}
= \overline{\R}_{{}_{G\!H}} \, $.  This raises the
question of whether  $ \overline{\R}_{{}_{W\!X}} $  and 
$ \overline{\R}_{{}_{G\!H}} $  do always coincide: this
problem is tackled and solved in the forthcoming paper
[EGH], as a byproduct of a more general uniqueness result
for braidings in  $ \fgstar $,  which on its own follows
from uniqueness of ``lifts'' of classical  $ \r $--matrices 
(a new notion which makes precise the idea of ``integrating''
a classical  $ \r $--matrix).   
                                             \par  
   The paper is organized as follows.  Section 1 is devoted to
recall some notions and results of quantum theory.  Section 2 
deals with the construction of braidings via quantum groups, after 
[GH]: in particular we point out its "compatibility" with the 
equivalence relation for quantisations, we prove the triviality of 
the associated infinitesimal braiding, and we sketch some 
examples.  Section 3 deals with the geometrical construction of 
braidings after [WX]: in particular we reformulate some results 
from  [{\it loc.~cit.}]  to make them fit with our language, and 
we prove that the associated infinitesimal braiding is trivial.  
Finally, section 4 is devoted to explicit computation of both  $ 
\overline{\R}_{{}_{W\!X}} $  and $ \overline{\R}_{{}_{G\!H}} $,  
which shows they do coincide.  

\vskip2,1truecm

\centerline {\bf \S\; 1. \ Definitions and recalls
from quantum group theory }  

\vskip13pt 

  {\bf 1.1 Topological  $ \kh $--modules  and topological
Hopf  $ \kh $--algebras.} \, Let  $ \Bbbk $  be a fixed field
of zero characteristic,  $ \hbar $  an indeterminate.  The
ring  $ \kh $  will always be considered as a topological
ring w.r.t.~the  $ \hbar $--adic  topology.  Let  $ X $ be
any  $ \kh $--module.  We set  $ \, X_0 := X \big/ \hbar X
= \Bbbk \otimes_{\kh} X \, $,  \, a  $ \Bbbk $--module  (via
scalar restriction $ \, \kh \rightarrow \kh \big/ \hbar \, \kh
\cong \Bbbk \, $)  which we call the  {\sl specialisation\/} 
of  $ X $  at  $ \, \hbar = 0 \, $,  or  {\sl semiclassical
limit\/}  of  $ X \, $;  we shall also use notation  $ \, X
\,{\buildrel \, \hbar \rightarrow 0 \, \over \llongrightarrow}\,
\overline{Y} \, $  to mean  $ \, X_0 \cong \overline{Y} \, $. 
For later use, we also set  $ \, \lsF X := \Bbbk((h))
\otimes_{\kh} X \, $,  \, a vector space over 
$ \Bbbk((h)) \, $.  
%  
%  
% If  $ X $  is a topological  $ \kh $--module,  we let its 
% {\it full dual\/}  to be  $ \, X^* := \text{\it Hom\,}_{\kh}
% \big( X, \kh \big) \, $,  \, and its  {\it topological dual\/} 
% to be  $ \, X^\star := \big\{\, f \in X^* \,\big\vert\, f \text{\
% is continuous} \,\big\} \, $.  Note that  $ \, X^* = X^\star \, $ 
% when the topology on  $ X $  is the  $ \hbar $--adic  one.  
%  
%  
%                                             \par
%    We introduce now two tensor categories of topological 
% $ \kh $--modules,  $ {\Cal T}_{\otimeshat} $  and 
% $ {\Cal P}_{\otimestilde} $:  the first one is modeled on the tensor
% category of  {\sl discrete\/}  topological  $ \Bbbk $--vector 
% spaces, the second one is modeled on the category of  {\sl
% linearly compact\/}  topological  $ \Bbbk $--vector  spaces.  
%   
%   
                                            \par
   Let  $ {\Cal T}_{\otimeshat} $  be the category whose objects
are all topological  $ \kh $--modules  which are topologically free
(i.e.~isomorphic to  $ V[[\hbar]] $  for some  $ \Bbbk $--vector 
space $ V $,  with the  $ \hbar $--adic  topology) and whose 
morphisms are the  $ \kh $--linear  maps (which are automatically 
continuous). This is a tensor category w.r.t.~the tensor product  
$ \, T_1 \otimeshat T_2 \, $  defined to be the separated 
$ \hbar $--adic  completion of the algebraic tensor product 
$ \, T_1 \otimes_{\kh} T_2 \, $ (for all  $ T_1 $,  $ T_2
\in {\Cal T}_{\otimeshat} $).  
                                            \par
   Let  $ {\Cal P}_{\otimestilde} $  be the category whose objects
are all topological  $ \kh $--modules  isomorphic to modules of 
the type  $ {\kh}^E $  (the Cartesian product indexed by $ E $,  
with the Tikhonov product topology) for some set  $ E \, $: \, 
these are complete w.r.t.~to the weak topology, in fact they are 
isomorphic to the projective limit of their finite free submodules 
(each one taken with the  $ \hbar $--adic  topology); the
morphisms in $ {\Cal P}_{\otimestilde} $  are the  $ \kh $--linear 
continuous maps.  This is a tensor category w.r.t.~the tensor 
product  $ \, P_1 \otimestilde P_2 \, $  defined to be the 
completion of the algebraic tensor product  $ \, P_1 \otimes_{\kh} 
P_2 \, $  w.r.t.~the weak topology:  therefore  $ \, P_i \cong 
{\kh}^{E_i} $  ($ i = 1 $, $ 2 $)  yields  $ \, \, P_1 
\otimestilde P_2 \cong {\kh}^{E_1 \times E_2} \, $  (for
all  $ P_1 $,  $ P_2 \in {\Cal P}_{\otimestilde} $).  
                                            \par  
   Note that the objects of  $ {\Cal T}_{\otimeshat} $  and of
$ {\Cal P}_{\otimestilde} $  are complete and separated w.r.t.~the 
$ \hbar $--adic  topology, whence one has  $ \, X \cong 
X_0[[\hbar]] \, $ (as  $ \kh $--modules)  for each of them.   
                                            \par  
   To simplify notation, in the sequel we shall usually write
simply  $ \, \otimes \, $  \hbox{for either  $ \, \otimeshat \, $ 
or  $ \, \otimestilde \, $.}   

\vskip9pt 

\proclaim{Definition 1.2}  (cf.~[Dr], \S~7) 
                                         \hfill\break
  \indent  (a) \, We call  {\sl quantized universal enveloping
algebra}  (in short, QUEA)  any Hopf algebra  $ \, H \, $  in the 
category  $ \, {\Cal T}_{\otimeshat} \, $  such that  $ \, H_0 := 
H \big/ \hbar H \, $  is a co-Poisson Hopf algebra isomorphic to  
$ U(\gerg) $  for some finite dimensional Lie bialgebra $ \gerg $  
(over  $ \Bbbk $);  in this case we write  $ \, H = U_\hbar(\gerg) \, 
$,  \, and say  $ H $  is a  {\sl quantisation} of  $ U(\gerg) $.  
We call \QUEA{} the subcategory of $ {\Cal T}_{\otimeshat} $  
whose objects are QUEA (relative to all possible  $ \gerg $),  
with the obvious morphisms.   
                                         \hfill\break
  \indent  (b) \, We call  {\sl quantized formal series Hopf algebra}
(in short, QFSHA)  any Hopf algebra  $ \, K \, $  in the category 
$ \, {\Cal P}_{\otimestilde} \, $  such that  $ \, K_0 := K \big/ 
\hbar K \, $  is a topological Poisson Hopf algebra isomorphic to 
$ F[[\gerg]] $  for some finite dimensional Lie bialgebra  $ \gerg 
$ (over  $ \Bbbk $);  then we write  $ \, H = F_\hbar[[\gerg]] \, $,  
\, and say  $ K $  is a  {\sl quantisation}  of  $ F[[\gerg]] $.  
We call \QFSHA{} the full subcategory of  $ {\Cal 
P}_{\otimestilde} $  whose objects are QFSHA (relative to all 
possible  $ \gerg $),  with the obvious morphisms.  
                                         \hfill\break
  \indent  (c) \, If  $ H_1 $,  $ H_2 $,  are two quantisations of
$ U(\gerg) $,  resp.~of  $ F[[\gerg]] $  (for some Lie bialgebra $ 
\gerg $),  we say that  {\sl  $ H_1 $  is equivalent to  $ H_2 $}, 
and we write  $ \, H_1 \equiv H_2 \, $,  if there is an 
isomorphism $ \, \varphi \, \colon H_1 \cong H_2 \, $  (in \QUEA, 
resp.~in \QFSHA) and a  $ \kh $--linear  isomorphism  $ \, 
\varphi_+ \, \colon H_1 \cong H_2 \, $  such that  $ \, \varphi = 
\hbox{\rm id} + \hbar \, \varphi_+ \, $  when identifying  $ H_1 $ 
and  $ H_2 $  with $ U(\gerg)[[\hbar]] $,  resp.~with  $ 
F[[\gerg]][[\hbar]] $. 
\endproclaim  

\vskip9pt 

  {\bf Remarks 1.3:} \, {\it (a)} \,  Note that the objects
of  \QUEA{}  and of  \QFSHA{}  are  {\sl topological\/}  Hopf 
algebras, not standard ones.  As a matter of notation, if  $ H $ 
is any Hopf algebra (maybe topological), we shall denote by  $ m $ 
its product, by  $ 1 $  its unit element, by  $ \Delta $  its 
coproduct, by  $ \epsilon $  its counit and by  $ S $  its 
antipode (with a subscript  $ H $  if necessary).  
                                                   \par
   {\it (b)} \, If  $ \, H \in \HA_{\otimeshat} $
is such that  $ \, H_0 := H \big/ \hbar H \, $  as a Hopf algebra 
is isomorphic to  $ U(\gerg) $  for some Lie algebra  $ \gerg $,  
then $ \, H_0 = U(\gerg) \, $  is also a  {\sl co-Poisson\/}  Hopf 
algebra, w.r.t.~the Poisson cobracket  $ \delta $  defined as
follows: if  $ \, x \in H_0 \, $  and  $ \, x' \in H \, $  gives 
$ \, x = x' + h \, H \, $,  \, then  $ \, \delta(x) := \big( h^{-1}
\, \big( \Delta(x') - \Delta^{\text{op}}(x') \big) \big) + h \,
H \otimeshat H \, $;  \, then (by [Dr], \S 3, Theorem 2) the
restriction of  $ \delta $  makes  $ \gerg $  into a Lie
bialgebra.  Similarly, if  $ \, K \in \HA_{\otimestilde} $  is
such that  $ \, K_0 := K \big/ \hbar K \, $  is a topological
Poisson Hopf algebra isomorphic to  $ F[[\gerg]] $  for some Lie
algebra  $ \gerg $  then  $ \, K_0 = F[[\gerg]] \, $  is also a
topological  {\sl Poisson\/}  Hopf algebra, w.r.t.~the Poisson
bracket  $ \{\,\ ,\ \} $  defined as follows: if  $ \, x $, 
$ y \in K_0 \, $  and  $ \, x' $,  $ y' \in K \, $  give  $ \,
x = x' + h \, K $,  $ \, y = y' + h \, K $,  \, then  $ \,
\{x,y\} := \big( h^{-1} (x' \, y' - y' \, x') \big) + h \,
K \, $;  \, then  $ \gerg $  is a bialgebra again. These
natural co-Poisson and Poisson structures are the ones
considered in Definition 1.2 above.  
                                                   \par
   {\it (c)} \, Clearly \QUEA, resp.~\QFSHA,  is a  {\sl
tensor\/}  subcategory of  $ {\Cal T}_{\otimeshat} $, 
resp.~of  $ {\Cal P}_{\otimestilde} $. 
                                                   \par
   {\it (d)} \, We make a finiteness assumption on  $ \, \text{\it
dim}\,(\gerg) \, $,  \, but infinite-dimensional cases can also be 
"reasonably" handled as explained in [Ga5], \S 3.9.   

 \vskip9pt 
%  
%  
%  \eject  
%  
%  

  {\bf 1.5 Drinfeld's functors.} \,  Let  $ H $  be a Hopf
algebra (of any type) over  $ \kh $.  For each  $ \, n \in \N $,  
define  $ \; \Delta^n \colon H \longrightarrow H^{\otimes n} \; $  
by  $ \, \Delta^0 := \epsilon \, $,  $ \, \Delta^1 := \hbox{\rm
id}_{\scriptscriptstyle H} $,  \, and  $ \, \Delta^n := \big( 
\Delta \otimes \hbox{\rm id}_{\scriptscriptstyle H}^{\,\otimes
(n-2)} \big) \circ \Delta^{n-1} \, $  if  $ \, n \geq 2 $.  Then
set  $ \; \delta_n = {(\hbox{\rm id}_{\scriptscriptstyle H} -
\epsilon)}^{\otimes n} \circ \Delta^n \, $,  \, for all  $ \,
n \in \N_+ \, $.  Finally, define   
  $$  H' := \big\{\, a \in H \,\big\vert\; \delta_n(a) \in
\hbar^n H^{\otimes n} \; \forall\, n \in \N \,\big\}  \qquad 
\big( \subseteq H \, \big ) \; .   $$  
   \indent   Now let  $ \, I_{\scriptscriptstyle H} :=
{\epsilon_{\scriptscriptstyle H}}^{\hskip-3pt -1} \big( \hbar \, 
\kh \big) \, $;  \, set  $ \, H^\times := \sum\limits_{n \geq 0} 
\hbar^{-n} {I_{\scriptscriptstyle H}}^{\!n} = \bigcup\limits_{n 
\geq 0} {\big( \hbar^{-1} I_{\scriptscriptstyle H} \big)}^n \, $  
$ \, \big( \subseteq \lsF H \, \big) $,  \, and 
  $$  H^\vee :=  \text{(separated)  $ \hbar $--adic  completion
of the  $ \kh $--module }  H^\times \; .  $$   

\vskip9pt 

   The following is the first important result we need:  

\vskip9pt 

\proclaim {Theorem 1.6} ("The quantum duality principle"; 
cf.~[Ga5], Theorem 1.6)   
                                        \hfill\break
   \indent   The assignments  $ \, H \mapsto H^\vee \, $  and
$ \, H \mapsto H' \, $  respectively define functors of tensor 
categories  $ \, \QFSHA \longrightarrow \QUEA \, $  and  $ \, 
\QUEA \longrightarrow \QFSHA \, $.  These functors are inverse to 
each other.  Indeed,  for all  $ \, U_\hbar(\gerg) \in \QUEA \, $ 
and all  $ \, F_\hbar[[\gerg]] \in \QFSHA \, $  one has 
  $$  {U_\hbar(\gerg)}' \Big/ \hbar \, {U_\hbar(\gerg)}' =
F[[\gerg^*]] \, , \qquad  {F_\hbar[[\gerg]]}^\vee \Big/ \hbar \, 
{F_\hbar[[\gerg]]}^\vee = U(\gerg^*)  $$   
(where  $ \gerg^* $  is the dual to  $ \gerg $),  i.e.~$ \,
{U_\hbar(\gerg)}' \! = \! F_\hbar[[\gerg^*]] \, $  and  $ \,
{F_\hbar[[\gerg]]}^\vee \! = U_\hbar(\gerg^*) \, $.  Moreover,
the functors preserve equivalence,  i.e.~$ \, H_1 \equiv H_2
\, $  implies  $ \, {H_1}^{\!\vee} \equiv {H_2}^{\!\vee} \, $ 
or  $ \, {H_1}' \equiv {H_2}' \, $.   \qed   
\endproclaim  

\vskip9pt  

   {\bf 1.7 An explicit description of  $ \, {U_\hbar(\gerg)}' \, $.} \,
Given any QUEA, say  $ U_\hbar(\gerg) $,  we can give a rather 
explicit description of  $ \, {U_\hbar(\gerg)}' \, $.  In fact, 
one has (see [Ga5], \S 3.5):  

\vskip3pt 

% 
%    {\it  Given any basis of  $ \gerg $,  there exists a lift  $ \,
% {\{x_i\}}_{i \in \Cal{I}} \, $  of it in  $ U_\hbar(\gerg) $  such that 
% $ \, \epsilon(x_i) = 0 \, $  and  $ {U_\hbar(\gerg)}' $  is nothing but
% the topological  $ \kh $--algebra in  $ \Cal{P}_{\otimestilde} $ 
% generated (in topological sense) by  $ \, {\{\hbar \, x_i\big\}}_{i
% \in \Cal{I}} \, $,  \, thus  $ \, {U_\hbar(\gerg)}' = \Big\{
% \sum_{\underline{e} \in {(\N^\Cal{I})}_0} \hskip-1pt
% a_{\underline{e}} \, \hbar^{|\underline{e}|} \, x^{\underline{e}}
% \hskip5pt \Big| \hskip3pt a_{\underline{e}} \in \kh
% \;\, \forall \; \underline{e} \,\Big\} \, $  as a
% subset of  $ \, U_\hbar(\gerg) $.}  

   {\it  Given a basis  $ \, \big\{ \overline{x}_1, \dots,
\overline{x}_d \big\} \, $  of  $ \gerg $,  there is a lift $ \, 
\big\{x_1, \dots, x_d\big\} \, $  of it in  $ U_\hbar(\gerg) $ 
such that  $ \, \epsilon(x_i) = 0 \, $  and  $ {U_\hbar(\gerg)}' $ 
is just the topological  $ \kh $--algebra  in $ 
\Cal{P}_{\otimestilde} $  generated (topologically) by  $ \, 
\big\{\hbar \, x_1, \dots, \hbar \, x_d\big\} \, $;  \, so  $ \, 
{U_\hbar(\gerg)}' = \Big\{ \sum_{\underline{e} \in \N^d}  
\hskip-1pt  a_{\underline{e}} \, \hbar^{|\underline{e}|} \, 
x^{\underline{e}}  \hskip5pt  \Big| \hskip3pt a_{\underline{e}} 
\in \kh \;\, \forall \; \underline{e} \,\Big\} \, $  as a 
\hbox{subset of  $ U_\hbar(\gerg) $.}}  

\vskip3pt 

   Hereafter, we use notation  $ \, x^{\,\underline{e}} :=
\prod\limits_{i=1}^d x_i^{\,\underline{e}_i} \, $  and  $ \, 
|\underline{e}| := \sum\limits_{i=1}^d e_i \, $  for  \hbox{all $ 
\, \underline{e} = \big( e_1, \dots, e_d \big) \in \N^d $.}  

\vskip9pt 

\proclaim{Definition 1.8}  (cf.~[Dr1], [CP], [Re])  
                                         \hfill\break
  \indent  (a) \, A Hopf algebra  $ H $  (in any tensor category) is
called  {\sl quasitriangular}  if there is  $ \, R \in H \otimes H 
\, $  (tensor product within the category), called  {\sl the $ R 
$--matrix  of}  $ H $,  such that  
  $$  \hbox{ $ \eqalign{ 
   R \cdot \Delta (a) \cdot R^{-1} =  &  \; {\hbox{\rm Ad}}(R)
(\Delta(a)) = \Delta^{\text{op}}(a)  \cr 
   (\Delta \otimes \hbox{\rm id}) (R) = R_{13} R_{23} \,  &  , 
\qquad  (\hbox{\rm id} \otimes \Delta) (R) = R_{13} R_{12} 
\cr } $ }  \eqno (1.1)  $$   
% 
%   $$ \; R \cdot \Delta (a) \cdot R^{-1} = {\hbox{\rm Ad}}(R)
% (\Delta(a)) = \Delta^{\text{op}}(a) \, , 
% \; (\Delta \otimes \hbox{\rm id}) (R) = R_{13} R_{23} \, , 
% \; (\hbox{\rm id} \otimes \Delta) (R) = R_{13} R_{12} \, ,  $$  
% 
where  $ \, \Delta^{\text{op}} := \sigma \circ \Delta(a) \, $  with 
$ \, \sigma \colon \, H^{\otimes 2} \to H^{\otimes 2} \, $,  $ \,
a \otimes b \mapsto b \otimes a  \, $,  \, and  $ \, R_{12}, R_{13},
R_{23} \in H^{\otimes 3} \, $,  $ \, R_{12} = R \otimes 1 \, $, 
$ \, R_{23} = 1 \otimes R \, $,  $ \, R_{13} = (\sigma \otimes
\hbox{\rm id}) (R_{23}) =  (\hbox{\rm id} \otimes \sigma) (R_{12})
\, $.  The algebra is called  {\sl triangular},  and the 
$ R $--matrix {\sl unitary},  if in addition  $ \, R^{-1}
= R^{\text{op}} := \sigma(R) \, $.   
                                                      \par   
   We call  \QTQUEA,  resp.~\TQUEA,  the subcategory of  \QUEA{} 
whose objects are all the quasitriangular, resp.~the triangular, 
QUEA (in short QTQUEA, resp.~TQUEA) and whose morphisms  $ \, 
\varphi \, \colon \, H_1 \longrightarrow H_2 \, $  enjoy $ \, 
\phi^{\otimes 2}(R_1) = R_2 \, $.  
                                               \par   
%                                          \hfill\break
%  
%  
%  \eject  
%  
%  
%   \indent  
% 
% 
   (b) \, A Hopf algebra  $ H $  (in any tensor category) is
called  {\sl braided}  if there is an  {\sl algebra automorphism} 
$ \; \R \, \colon \, H \otimes H \longrightarrow H \otimes H \; $ 
in the category, called  {\sl the braiding operator  {\text (or 
simply}  the braiding{\text )}  of}  $ H $,  different from  $ \; 
\sigma \colon a \otimes b \mapsto b \otimes a \; $  and such that  
  $$  \hbox{ $ \eqalign{ 
   \R \, \circ  &  \, \Delta = \Delta^{\text{op}}  \cr  
   (\Delta \otimes \hbox{\rm id}) \circ \R = \R_{13}
\circ \R_{23} \circ (\Delta \otimes \hbox{\rm id}) \, ,  &  \qquad 
   (\hbox{\rm id} \otimes \Delta) \circ \R = \R_{13}
\circ \R_{12} \circ (\hbox{\rm id} \otimes \Delta)  \cr } $ }  
\eqno (1.2)  $$    
where  $ \R_{12}, \R_{13}, \R_{23} $  are the automorphisms of 
$ H^{\otimes 3} $  defined by  $ \, \R_{12} = \R \otimes \hbox{\rm
id} \, $,  $ \, \R_{23} = \hbox{\rm id} \otimes \R \, $,  $ \,
\R_{13} = (\sigma \otimes \hbox{\rm id}) \circ (\hbox{\rm id}
\otimes \R) \circ (\sigma \otimes \hbox{\rm id}) \, $.  Moreover,
the braiding operator is said to be  {\sl unitary}  and the algebra
to be  {\sl rigid}  if in addition  $ \, \R^{-1} = \sigma \circ \R
\circ \sigma \, $.   
                                                      \par   
   We call  \BQFSHA,  resp.~\RBQFSHA,  the subcategory of  \QFSHA{} 
whose objects are all the braided, resp.~the rigid braided, QFSHA 
(in short BQFSHA, resp.~RBQFSHA) and whose morphisms  $ \, \psi \, 
\colon \, H_1 \longrightarrow H_2 \, $  enjoy  $ \, \psi^{\otimes 
2} \circ \R_1 = \R_2 \circ \psi^{\otimes 2} \, $.  
                                         \hfill\break
  \indent  (c) \, Let  $ \, \big(H_1,R_1\big) $,  $ \big(H_2,R_2\big)
\in \QTQUEA $.  We say that  {\sl  $ \, \big(H_1,R_1\big) \, $  is 
equivalent to  $ \, \big(H_2,R_2\big) \, $},  \, and we write  $ 
\, \big(H_1,R_1\big) \equiv \big(H_2,R_2\big) \, $,  \, if  $ \, 
H_1 \equiv H_2 \, $  in  \QUEA{}  via an equivalence  $ \, \varphi 
\, \colon H_1 \cong H_2 \, $  which is also an isomorphism in 
\QTQUEA{}  (i.e.~such that  $ \, \phi^{\otimes 2}(R_1) = R_2
\, $).  
                                         \hfill\break
  \indent  (d) \, Let  $ \, \big(H_1,\R_1\big) $,  $ \big( H_2,
\R_2 \big) \in \BQFSHA $.  We say that  {\sl  $ \, \big( H_1, \R_1
\big) \, $  is equivalent to  $ \, \big(H_2,\R_2\big) \, $},  \,
and we write  $ \, \big(H_1,\R_1\big) \equiv \big(H_2,\R_2\big)
\, $,  \, if  $ \, H_1 \equiv H_2 \, $  in  \QUEA{}  via an
equivalence which is also an isomorphism in  \BQFSHA{} 
(i.e.~such that  $ \, \psi^{\otimes 2} \circ \R_1 =
\R_2 \circ \psi^{\otimes 2} \, $).   
\endproclaim  

\vskip9pt 

   {\bf Remarks 1.9.} \, {\it (a)} \, It follows immediately from
(1.1) that  $ R $  {\it is a solution of the quantum Yang-Baxter 
equation\/  {\rm (in short, QYBE)}  in}  $ H^{\otimes 3} $, 
namely $ \; R_{12} R_{13} R_{23} = \! R_{23} R_{13} R_{12} \, $. 
This is the starting point for defining a braid group action on 
the tensor products of  $ H $--modules,  and then for constructing 
link invariants, following [Tu] (see also [CP], \S 15).  
                                          \par   
   Similarly, it follows from (1.2) that  $ \R $  {\it is a
solution of the QYBE in}  $ \text{\it End}\,(H^{\otimes 3}) $, 
namely  $ \, \R_{12} \circ \R_{13} \circ \R_{23} = \R_{23} \circ
\R_{13} \circ \R_{12} \, $.  Again, this implies the existence of
a braid group action on the tensor powers of  $ H $,  from which
one can start a search for link invariants.  
                                          \par   
   {\it (b)} \, It is proved in [EK] that, for any Lie bialgebra 
$ \gerg $,  there exists a QUEA, which we'll denote  $ \uhg $, 
whose semiclassical limit is isomorphic to  $ \ug $;  moreover, 
one has an identification  $ \, \uhg \cong \ug [[\hbar]] \, $  as 
$ \kh $--modules,  hence also  $ \, \uhg \otimes \uhg \cong \big( 
\ug \otimes \ug \big) [[\hbar]] \, $.  Here, like elsewhere in the 
following, the tensor products among  $ \kh $--modules  are 
topological tensor products.  In addition, if  $ \gerg $  is 
quasitriangular   --- as a Lie bialgebra ({\it cf.}~[CP]) ---   
and $ \r $  is its  $ \r $--matrix,  then there exists such a  $ 
\uhg $  which is quasitriangular as well   --- as a Hopf algebra 
---   with an  $ R $--matrix  $ \, R_\hbar \,\, (\, \in \uhg \otimes
\uhg \, $)  such that  $ \; R_\hbar \equiv 1 + \hbar \, \r \; 
\mod\, \hbar^2 \, $,  that is to say  $ \, R_\hbar = 1 + \hbar
\, \r + \O \left( \hbar^2 \right) \, $  with  $ \, \O \left( 
\hbar^2 \right) \in \hbar^2 \cdot \uhg \otimes \uhg \, $.  

\vskip2,1truecm

\centerline {\bf \S\; 2. \  Braidings from deformation 
quantisation } 

\vskip13pt 

\proclaim{Theorem 2.1}  ([GH], Th\'eor\`eme 2.1) Let  $ H $  be a 
QTQUEA, and let  $ R $  be its  $ R $--matrix.  Then the inner 
automorphism  $ \; {\hbox{\rm Ad}} (R) \, \colon \, H \otimes H 
\llongrightarrow H \otimes H \; $  of  $ \, H \otimes H \, $ 
restricts to an automorphism of  $ \, H' \otimes H' $,  and the 
pair  $ \, \Big( H', \, {\hbox{\rm Ad}}(R){\big\vert}_{H' \otimes 
H'} \Big) \, $  is a BQFSHA.   \qed  
\endproclaim 

\vskip9pt 

   As a first goal in this section we provide some further details
about Theorem 2.1:   

\vskip9pt 

\proclaim{Theorem 2.2}   
                                          \hfill\break   
   \indent   (a) \, The functor  $ \; {(\ )}' \, \colon \, \QUEA
\llongrightarrow \QFSHA \; $  yields by restriction two functors  
                                          \hfill\break   
  $$  {(\ )}' \, \colon \, \QTQUEA \llongrightarrow \BQFSHA \;\; , 
\qquad  \big(H,R\big) \mapsto \Big( H', \, {\hbox{\rm Ad}}(R)
{\big\vert}_{H' \otimes H'} \Big) \; \phantom{.}  $$   
  $$  {(\ )}' \, \colon \, \TQUEA \llongrightarrow \RBQFSHA \;\; , 
\qquad  \big(H,R\big) \mapsto \Big( H', \, {\hbox{\rm Ad}}(R)
{\big\vert}_{H' \otimes H'} \Big) \; .  $$   
   \indent   (b) \, The functors in (a) preserves equivalence
classes, i.e.~if  $ \, \big(H_1,R_1\big) \equiv \big(H_1,R_1\big) 
\, $  in  \QTQUEA{}  then  $ \, \Big( H', \, {\hbox{\rm Ad}}
(R){\big\vert}_{H' \otimes H'} \Big) \equiv \Big( H', \, 
{\hbox{\rm Ad}}(R){\big\vert}_{H' \otimes H'} \Big) \, $ 
in  \BQFSHA.  
\endproclaim  

\demo{Proof}  {\it (a)} \, Theorem 2.1 tells that the functor $ \; 
{(\ )}' \, \colon \, \QTQUEA \loongrightarrow \BQFSHA \; $ is 
well-defined on objects.  Moreover, if  $ \, \phi \, \colon \, 
\big(H_1,R_1\big) \longrightarrow \big(H_2,R_2\big) \, $  is a 
morphism in  \QTQUEA{}  then  $ \, \phi^{\otimes 2}(R_1) = R_2 
\, $,  \, whence  $ \, \phi^{\otimes 2} \circ {\hbox{\rm Ad}}(R_1)
{\big\vert}_{{(H'_1)}^{\otimes 2}} = {\hbox{\rm Ad}}(R_2)
{\big\vert}_{{(H'_2)}^{\otimes 2}} \circ \phi^{\otimes 2} \, $ 
follows at once, hence  $ \; \phi' := \phi{\big\vert}_{H'_1} \, 
\colon \, H'_1 \loongrightarrow H'_2 \; $  is a morphism in  
\BQFSHA. In addition, if  $ \, \big(H,R\big) \in \TQUEA \, $ 
then  $ \, R^{-1} = \sigma(R) \, $  yields  $ \, {\Big(
{\hbox{\rm Ad}}(R){\big\vert}_{{(H')}^{\otimes 2}} \Big)}^{-1} = 
{\hbox{\rm Ad}}\big( R^{-1} \big){\Big\vert}_{{(H')}^{\otimes 2}} 
= {\hbox{\rm Ad}}\big( \sigma(R) \big){\Big\vert}_{{(H')}^{\otimes 
2}} = \sigma \circ {\hbox{\rm Ad}}\big( \sigma(R) \big)%  
{\Big\vert}_{{(H')}^{\otimes 2}} \circ \sigma \, $, 
\, hence  $ \, {\hbox{\rm Ad}}\big( \sigma(R) \big)%  
{\Big\vert}_{{(H')}^{\otimes 2}} \, $  is unitary, q.e.d. 
                                                \par   
   {\it (b)} \, This follows easily from  {\it (a)\/}  and
the very definitions.   \qed   
\enddemo  

\vskip9pt 

   Second, as a consequence of Theorem 2.1 along with
the existence of quasitriangular quantisation of any
quasitriangular Lie bialgebra (cf.~[EK]) one gets
a braiding on  $ \, F[[\gerg^*]] \, $:  

\vskip9pt 

\proclaim{Corollary 2.3}  ([GH], Th\'eor\`eme 2.2)  \, Let 
$ \gerg $ be a (finite dimensional) quasitriangular Lie
bialgebra.  Then the topological Poisson Hopf algebra 
$ \, F[[\gerg^*]] \, $  is braided (in particular, its
braiding is a Poisson automorphism).  Moreover, there is
a quantisation of  $ \, F[[\gerg^*]] \, $  which is a braided
Hopf algebra whose braiding operator specialises into that of 
$ \, F[[\gerg^*]] \, $.   \qed   
\endproclaim  

\vskip9pt 

   {\bf 2.4 The triviality of the infinitesimal braiding.} \, Let 
$ \gerg $  and  $ \gerg^* $  be finite dimensional Lie bialgebras 
dual to each other.  Assume  $ \fgstar $  is braided (as a Poisson 
Hopf algebra),  $ \R $  being its braiding (which is a Poisson 
automorphism also).  Let  $ \, \frak{m}^\otimes_e \, $  be the 
(unique) maximal ideal of  $ \, F[[\gerg^* \oplus \gerg^*]] = 
F[[\gerg^*]] \otimes F[[\gerg^*]] \, $  (topological tensor 
product, after [Di], Ch.~1).  Since  $ \R $  is an algebra 
automorphism,  $ \, \R \big( \frak{m}^\otimes_e \big) = 
\frak{m}^\otimes_e \, $  and  $ \R $  induces an automorphism $ 
\overline{\R} $  of the vector space  $  \frak{m}^\otimes_e \Big/ 
{\big( \frak{m}^\otimes_e \big)}^{\!2} \, $.  Now,  $ \, 
\frak{m}^\otimes_e \Big/ {\big( \frak{m}^\otimes_e \big)}^{\!2} \, 
$  with the Lie bracket induced by the Poisson bracket of $ \, 
F[[\gerg^* \oplus \gerg^*]] \, $  identifies with the Lie algebra  
$ \, \gerg \oplus \gerg \, $;  \, since  $ \R $  is also an 
automorphism of Poisson algebras, the map  $ \overline{\R} $  is 
an automorphism of the Lie algebra  $ \, \gerg \oplus \gerg \, $; 
\, of course  $ \overline{\R} $  inherits also other properties of 
the braiding  $ \R $,  in particular  $ \R $  and  $ \overline\R $ 
are solutions of the QYBE, hence we call it  {\it the 
infinitesimal braiding\/}  associated to  $ \R $.  
                                                \par   
   Now assume in addition that  $ \gerg $  be quasitriangular, and
the braiding  $ \R $  on  $ \fgstar $  is provided as in Corollary 
2.3.  Namely, let  $ \, \big( \uhg, R \big) \in \QTQUEA \, $  be a 
quantisation of the quasitriangular Lie bialgebra  $ \, (\gerg,\r) 
\, $:  by definition, this means that  $ \uhg $  has semiclassical 
limit (i.e.~specialisation at  $ \, \hbar = 0 $)  the co-Poisson 
Hopf algebra  $ \ug $  and, in the identification  $ \, \uhg = 
\ug[[\hbar]] \, $  (as topological  $ \kh $--modules),  $ \; R = 1 
+ \hbar \, \r + \O(\hbar^2) \, $  for some  $ \, \O(\hbar^2) \in 
\hbar^2 \, \ug[[\hbar]] \, $.  Then $ \R $  is the braiding of  $ 
\, F[[\gerg^* \oplus \gerg^*]] = {\uhg}' \otimes {\uhg}' \mod\, 
\hbar \, $  which is obtained as specialisation at  $ \, \hbar = 0 
\, $  of  $ \, {\hbox{\rm Ad}}(R){\Big\vert}_{{\uhg}' \otimes 
{\uhg}'} \, $,  \, thanks to Theorem 2.1.  Then our next result is 
that the associated infinitesimal braiding  $ \overline{\R} $  is 
always trivial:  

\vskip9pt 

\proclaim{Theorem 2.5} \, The infinitesimal braiding  $ \; 
\overline\R \, \colon \gerg \oplus \gerg \loongrightarrow \gerg 
\oplus \gerg \; $ is trivial, i.e.~$ \, \overline\R = \hbox{\rm
id}_{\gerg \oplus \gerg} \, $.  
\endproclaim  

\demo{Proof}  Let  $ \, \big\{\overline{x}_1, \dots, 
\overline{x}_d \big\} \, $  be a basis of  $ \gerg $,  and pick a 
lift  $ \, \big\{ x_1, \dots, x_d\big\} \, $  of it in  $ 
U_\hbar(\gerg) $  as explained in \S 1.7, so that  $ \, 
{U_\hbar(\gerg)}' = \Big\{ \sum_{\underline{e} \in \N^d}  
\hskip-1pt  a_{\underline{e}} \, \hbar^{|\underline{e}|} \, 
x^{\underline{e}} \, = \sum_{\underline{e} \in \N^d}  \hskip-1pt 
a_{\underline{e}} \, \tilde{x}^{\underline{e}}  \hskip5pt  \Big| 
\hskip3pt a_{\underline{e}} \in \kh \;\, \forall \; \underline{e} 
\,\Big\} \, $,  \; where  $ \, \tilde{x}_i := \hbar \, x_i \, $  
(for all  $ i \, $)  are topological generators of  $ 
{U_\hbar(\gerg)}' $. Then  $ \, {U_\hbar(\gerg)}' \otimes 
{U_\hbar(\gerg)}' \, $  is generated by the  $ \, {}_1\tilde{x}_i 
:= \tilde{x}_i \otimes 1 \, $  and the  $ {}_2\tilde{x}_i := 1 
\otimes \tilde{x}_i \, $,  \, for all  $ i $.  On the other hand, 
one has  $ \, U_\hbar(\gerg) = \Bbbk[x_1, \dots, x_d][[\hbar]] \, 
$  as topological  $ \kh $--modules, whence  $ \, U_\hbar(\gerg) 
\otimes U_\hbar(\gerg) = \big( \Bbbk[{}_1x_1, \dots, 
{}_1x_d,{}_2x_1, \dots, {}_2x_d] \big)[[\hbar]] \, $.  Then we 
have an $ \hbar $--adic  expansion of  $ R $  and of  $ R^{-1} $,  
namely  $ \; R = \sum\nolimits_{n\geq0} 
P^+_n\big({}_1\underline{x} \, ; {}_2\underline{x} \, \big) \, 
\hbar^n \, $,  $ \; R^{-1} = \sum\nolimits_{m \geq 0} 
P^-_m\big({}_1\underline{x} \, ; {}_2\underline{x} \, \big) \, 
\hbar^m \; $  for some polynomials $ P^+_n\!\big({}_1\underline{x} 
\, ; {}_2\underline{x} \, \big) \hskip-2,5pt  =  \hskip-2,5pt 
P^+_n\!\big({}_1x_1,\dots,{}_1x_d;{}_2x_1,\dots,{}_2x_d \big) $, $ 
P^-_m\big({}_1\underline{x} \, ; {}_2\underline{x} \, \big) 
\hskip-2,5pt  =  \hskip-2,5pt 
P^-_m\big({}_1x_1,\dots,{}_1x_d;{}_2x_1,\dots,{}_2x_d\big) $. Now, 
the condition  $ \, R = 1^\otimes + \hbar \, \r + 
\O\big(\hbar^2\big) \, $  (with  $ \, 1^\otimes := 1 \otimes 1 \, 
$)  forces  $ \; P_0^+ = 1 = P_0^- \, $,  $ \; P_1^+ = \sum_{i,j} 
c_{i,j} \cdot {}_1x_i \, {}_2x_j = - P_1^+ \; $  for some  $ \, 
c_{i,j} \in \Bbbk \, $  such that  $ \, \r = \sum_{i,j} c_{i,j}
\cdot \overline{x}_i \otimes \overline{x}_j \, $.  In addition,
any  $ R $--matrix  enjoys $ \, (\epsilon \otimes \hbox{\rm id}\,)
(R) = 1 = (\hbox{\rm id} \otimes \epsilon)(R) \, $,  \, hence also 
$ \, (\epsilon \otimes \hbox{\rm id}\,)\big(R^{-1}\big) = 1 = 
(\hbox{\rm id} \otimes \epsilon)\big(R^{-1}\big) \, $;  \, setting 
$ \, P_\pm := R^{\pm 1} - 1 \, $,  \, this implies  $ \; (\epsilon
\otimes \hbox{\rm id}\,)(P_\pm) = 0 = (\hbox{\rm id} \otimes
\epsilon)(P_\pm) \, $.   
                                                     \par   
   Now, for all  $ \ell $  consider  $ \; \big({\hbox{\rm Ad}}%  
(R)\big)({}_s\tilde{x}_\ell) = R \cdot {}_1\tilde{x}_i \cdot 
R^{-1} \, $:  \, we have  
  $$  \hbox{ $ \eqalign{ 
   \big({\hbox{\rm Ad}}(R)\big)({}_1\tilde{x}_\ell) =
R \cdot {}_1\tilde{x}_\ell \cdot R^{-1} = ( 1 + P_+ )  &  \cdot 
{}_1\tilde{x}_\ell \cdot ( 1 + P_- ) =   \hfill  \cr   
   \hfill   = {}_1\tilde{x}_\ell + P_+  &  \cdot {}_1\tilde{x}_\ell
+ {}_1\tilde{x}_\ell \cdot P_- + P_+ \cdot {}_1\tilde{x}_\ell 
\cdot P_- \; .  \cr } $ }   \eqno (2.1)  $$   We know that this 
element belongs to  $ \, {\big( U_\hbar(\gerg) \otimes 
U_\hbar(\gerg) \big)}' = \Bbbk[[{}_1\tilde{x}_1, \dots, 
{}_1\tilde{x}_d, {}_2\tilde{x}_1, \dots, {}_2\tilde{x}_d,h]]
\, $,  \, so we can write it as a series; since  $ \; (\epsilon
\otimes \hbox{\rm id}\,)(P_\pm) = 0 \, $  and  $ \, (\epsilon
\otimes \hbox{\rm id}\,) \, $  is a morphism we have  $ \;
(\epsilon \otimes \hbox{\rm id}\,) \big( P_+ \cdot
{}_1\tilde{x}_\ell + {}_1\tilde{x}_\ell \cdot P_- + P_+ \cdot
{}_1\tilde{x}_\ell \cdot P_- \big) = 0 \, $:  \, recalling that 
$ \, \epsilon\big({}_s\tilde{x}_i\big) = 0 \, $  this means that  
  $$  P_+ \cdot {}_1\tilde{x}_\ell + {}_1\tilde{x}_\ell \cdot P_-
+ P_+ \cdot {}_1\tilde{x}_\ell \cdot P_- = 
\sum\nolimits_{\underline{e}^{(1)},\underline{e}^{(2)} \in \N^d} 
\hskip-1pt  a_{\underline{e}^{(1)},\underline{e}^{(2)}} \, 
{}_1\tilde{x}^{\underline{e}^{(1)}} \, 
{}_2\tilde{x}^{\underline{e}^{(2)}}  $$   (where  $ \, 
a_{\underline{e}^{(1)},\underline{e}^{(2)}} \in \kh \, $  for all 
$ \underline{e}^{(1)} $,  $ \underline{e}^{(2)} \, $)  with  $ \,
a_{\underline{e}^{(1)},\underline{0}} = 0 =
a_{\underline{0},\underline{e}^{(2)}} \, $  for all 
$ \underline{e}^{(1)} $,  $ \underline{e}^{(2)} \, $, \, thus   
  $$  P_+ \cdot {}_1\tilde{x}_\ell + {}_1\tilde{x}_\ell \cdot P_-
+ P_+ \cdot {}_1\tilde{x}_\ell \cdot P_- = 
\sum\nolimits_{|\underline{e}^{(1)}|,|\underline{e}^{(2)}| > 1} 
\hskip-1pt  a_{\underline{e}^{(1)},\underline{e}^{(2)}} \, 
{}_1\tilde{x}^{\underline{e}^{(1)}} \, 
{}_2\tilde{x}^{\underline{e}^{(2)}} \; .   \eqno (2.2)  $$  Now  $ 
\, {}_1\tilde{x}^{\underline{e}^{(1)}} $, $ 
{}_2\tilde{x}^{\underline{e}^{(2)}} \! \mod \hbar \, {\big( 
{U_\hbar(\gerg)}^{\otimes 2} \big)}' \, $  belong to  $ \, 
\frak{m}^\otimes_e \, $  (notation of \S 2.4) as soon as $ \, 
|\underline{e}^{(1)}| > 1 \, $,  $ \, |\underline{e}^{(2)}| 
> 1 \, $;  \, so (2.2) gives  
  $$  \big( P_+ \cdot {}_1\tilde{x}_\ell + {}_1\tilde{x}_\ell
\cdot P_- + P_+ \cdot {}_1\tilde{x}_\ell \cdot P_- \big) \,
\equiv \, 0 \mod \hbar \, {\big( U_\hbar(\gerg)^{\otimes 2}
\big)}' \in {\big( \frak{m}^\otimes_e \big)}^2  $$   
and this along with (2.1) yields (for all  $ \, \ell = 1 $, 
$ \dots $,  $ d \, $)   
  $$  \displaylines{ 
   \R \, \bigg( \Big( {}_1\tilde{x}_\ell \! \mod \hbar \,
{\big( U_\hbar(\gerg)^{\otimes 2} \big)}' \Big) \! \mod {\big( 
\frak{m}^\otimes_e \big)}^2 \bigg) =   \hfill  \cr   
   = \Big( \big({\hbox{\rm Ad}}(R)\big)({}_1\tilde{x}_\ell) \!
\mod \hbar \, {\big( U_\hbar(\gerg)^{\otimes 2} \big)}' \Big) \! 
\mod {\big( \frak{m}^\otimes_e \big)}^2 =  \cr   
   \hfill   = \Big( {}_1\tilde{x}_\ell \! \mod \hbar \,
{\big( U_\hbar(\gerg)^{\otimes 2} \big)}' \Big) \! \mod {\big( 
\frak{m}^\otimes_e \big)}^2 \; .  \cr }  $$   Similarly one gets 
(for all  $ \, \ell = 1 $,  $ \dots $,  $ d \, $)  
  $$  \R \, \bigg( \Big( {}_2\tilde{x}_\ell \! \mod \hbar \,
{\big( U_\hbar(\gerg)^{\otimes 2} \big)}' \Big) \! \mod {\big( 
\frak{m}^\otimes_e \big)}^2 \bigg) = \Big( {}_2\tilde{x}_\ell \! 
\mod \hbar \, {\big( U_\hbar(\gerg)^{\otimes 2} \big)}' \Big) \! 
\mod {\big( \frak{m}^\otimes_e \big)}^2 \; .  $$  Letting  $ \; 
{}_s\check{x}_\ell := \Big( {}_s\tilde{x}_\ell \! \mod \hbar \, 
{\big( U_\hbar(\gerg)^{\otimes 2} \big)}' \Big) \! \mod {\big( 
\frak{m}^\otimes_e \big)}^2 \in \frak{m}^\otimes_e \Big/ {\big( 
\frak{m}^\otimes_e \big)}^2 = \gerg \oplus \gerg \; $  (for all  $ 
\, s=1 $,  $ 2 $  and  $ \, \ell = 1 $,  $ \dots $,  $ d \, $), \; 
we have in short  $ \; \R \big( {}_s\check{x}_\ell \big) = 
{}_s\check{x}_\ell \; $  for all  $ s $,  $ \ell $.  Since the $ 
{}_s\tilde{x}_\ell $  generate  $ {\big( U_\hbar(\gerg)^{\otimes 
2} \big)}' $,  the  $ {}_s\check{x}_\ell $  span  $ \, \gerg 
\oplus \gerg \, $,  \, hence we can conclude that  $ \R $  is 
trivial, as claimed.   \qed   
\enddemo  

\vskip9pt 

  {\bf 2.6 The example of semisimple and (untwisted) affine cases.}
\, In [Re] and [Ga1-2] the adjoint action of the  $ R $--matrix  
of the Jimbo-Lusztig's quantum groups  $ U_q(\gerg) $  was 
studied.  In this section we briefly outline how the results 
therein can be read as special occurrences of the ones cited here, 
namely the existence of braidings on  $ \uhg $.    
                                                \par   
   Let  $ \, \gerg = \gerg^\tau \, $  be a semisimple Lie algebra,
i.e.~a finite type Kac-Moody algebra, endowed with the Lie 
cobracket   --- depending on the parameter  $ \tau $  ---  
given in [Ga3], \S 1.3, which makes it into a Lie bialgebra;
in the following we shall also retain from  [{\it loc.~cit.}] 
all the notation we need: in particular, we denote by  $ Q $, 
resp.~$ P $,  the root lattice, resp.~the weight lattice, of 
$ \gerg $,  and by  $ r $  the rank of  $ \gerg $.  In particular,
when  $ \, \tau = 0 \, $  we have the standard Sklyanin-Drinfeld
cobracket. Similarly,  $ \gerg $  may be any untwisted affine
Kac-Moody algebra, as in [Ga4] (with corresponding notation).   
                                                \par   
   Now set  $ \, q := \exp(h) \, $;  then  $ \, \Bbbk(q) \, $ 
is a subring of  $ \kh $,  hence also all its subrings are. Let 
$ \, U_q(\gerg) \, $  be the Jimbo-Lusztig's quantum group over 
$ \Bbbk(q) $,  defined as  $ \, U_q(\gerg) := U_{q,
\varphi}^{\scriptscriptstyle Q}(\gerg) \, $  as in [Ga3],
\S 3.3, if  $ \gerg $  is finite, and as  $ \, U_q(\gerg) := 
U_q^{\scriptscriptstyle Q}(\gerg) \, $  as in [Ga4], \S 3.3, if 
$ \gerg $  is affine.  Furthermore, let  $ \, \uhat_q(\gerg) \, $ 
be the integer form of  $ U_q(\gerg) $  defined as  $ \, 
\uhat_q(\gerg) := {\frak U}_\varphi^{\scriptscriptstyle Q}(\gerg)
\, $  (over  $ \, A := \kqqm \, $)  as in [Ga3], \S 3.3, if 
$ \gerg $  is finite, and as  $ \, \uhat_q(\gerg) := {\frak
U}^{\scriptscriptstyle Q}(\gerg) \, $  (over the ring  $ A $  of 
rational functions in  $ q $  having no poles at roots of unity of 
odd order) as in [Ga4], \S 3.3, if  $ \gerg $  is affine.  In both
cases  $ \, A \, $  is a subring of  $ \Bbbk(q) $,  hence of 
$ \kh $,  thus we can define   
  $$  U_\hbar(\gerg) := \, \hbox{(separated)  $ \hbar $--adic 
completion of } \; \kh \otimes_A \uhat_q(\gerg) \; . 
\eqno (2.3)  $$   
   \indent   It is well known that  $ \, \uhat_q(\gerg) \Big/ (q-1) 
\, \uhat_q(\gerg) \cong U(\gerg) \, $:  this and (2.3) imply that 
$ \, U_\hbar(\gerg) \, $  has semiclassical limit  $ U(\gerg) $, 
thus it is a QUEA.  In fact,  $ U_\hbar(\gerg) $  is the well
known Drinfeld quantum group over  $ \kh $,  as defined in [Dr1],
\S 6.  In addition, let also  $ \, \utilde_q(\gerg) \, $  be
the integer form of  $ U_q(\gerg) $  defined as  $ \, \utilde_q
(\gerg) := {\Cal U}_\varphi^{\scriptscriptstyle Q}(\gerg) \, $ 
(over  $ \, A := \kqqm \, $)  as in [Ga3], \S 3.3, if 
$ \gerg $  is finite, and as $ \, \utilde_q(\gerg) :=
{\Cal U}^{\scriptscriptstyle Q}(\gerg) \, $  (over the
ring  $ A $  above)  as in [Ga4], \S 3.3, if  $ \gerg $ 
is affine.   
                                                \par   
   Similarly, we do the same for the dual Lie bialgebra 
$ \gerg^* $  (denoted  $ {\frak h} $  in  [{\it loc.~cit.}]), 
following [Ga3], \S 6   --- in the finite case ---   or [Ga4],
\S 5   --- in the affine case, thus getting  $ \, U_q(\gerg^*) $, 
$ \uhat_q(\gerg^*) $, $ \utilde_q(\gerg^*) $,  and  $ U_\hbar
(\gerg^*) $,  the last one being a QUEA with  $ U(\gerg^*) $ 
as semiclassical limit.  From the description in [Ga3--4], one 
sees that these objects are quite similar to the corresponding 
ones related to  $ \gerg $.  
                                                \par   
   Now consider  $ \; {\uhat_q(\gerg)}^{*} := \text{Hom}_A
\Big( \uhat_q(\gerg), A \Big) \, $; from [Ga3--4] we have the
identification  $ \; {\uhat_q}^{\,\, *}(\gerg) \cong \utilde_q
(\gerg^*) \, $, and also  $ \, \utilde_q(\gerg^*) \,{\buildrel
q \rightarrow 1 \over \llongrightarrow}\, \utilde_q(\gerg^*)
\Big/ (q-1) \, \utilde_q (\gerg^*) \cong F[[\gerg]] \, $. 
Thus letting   
  $$  F_\hbar[[\gerg]] := \, \hbox{(separated)  $ \hbar $--adic 
completion of }  \; \kh \otimes_A \utilde_q(\gerg^*)  
\eqno (2.4)  $$   
we have that  $ \, F_\hbar[[\gerg]] \, $  is a 
QFSHA, with semiclassical limit  $ F[[\gerg]] $.
                                                \par   
   The natural Hopf pairing  $ \; \langle\ ,\ \rangle \, \colon \,
\utilde_q(\gerg^*) \times \uhat_q(\gerg) \llongrightarrow A \; $ 
yields a Hopf pairing  $ \; \langle\ ,\ \rangle \, \colon \,
F_\hbar[[\gerg]] \times U_\hbar(\gerg) \llongrightarrow \kh \, $;  
moreover, it extends similarly to a perfect pairing  $ \; \langle
\ ,\ \rangle \, \colon \, U_q(\gerg^*) \times U_q(\gerg) 
\llongrightarrow \Bbbk(q) \; $.  The analysis in [Ga3--4]
shows that  
  $$  \utilde_q(\gerg) = {\big( \uhat_q(\gerg^*) \big)}^\circ
:= \Big\{\, y \in U_q(\gerg) \;\Big\vert\; \big\langle 
\uhat_q(\gerg^*), y \big\rangle \subseteq A \Big\} \, .   
\eqno (2.5)  $$   
In addition, by Proposition 1.4 we have also  
  $$  {U_\hbar(\gerg)}' = {\Big( \! {F_\hbar[[\gerg]]}^\vee
\Big)}^\circ := \Big\{\; y \in {}_{{}_F}U_\hbar(\gerg) 
\;\Big\vert\; \Big\langle {F_\hbar[[\gerg]]}^\vee, y \Big\rangle 
\subseteq \kh \;\Big\}  $$   
where we consider  $ \; \langle\ ,\ \rangle \, \colon \,
{}_{{}_F}F_\hbar[[\gerg]] \times {}_{{}_F}U_\hbar(\gerg)
\llongrightarrow \Bbbk((h)) \; $  to be the obvious pairing
obtained by scalar extension from  $ \, \langle\ , \ \rangle
\, \colon \, F_\hbar[[\gerg]] \times U_\hbar(\gerg)
\llongrightarrow \kh \, $.   
                                                \par   
   Now, the very definitions of all the objects involved yield
(via the analysis in [Ga3--4])  
  $$  \displaylines{ 
   {} \quad   {F_\hbar[[\gerg]]}^\vee = {\Big( \hbox{(separated) 
$ \hbar $--adic  completion of} \;\; \kh \otimes_A
\utilde_q(\gerg^*) \Big)}^\vee =   \hfill  \cr  
   {} \hfill   = \hbox{(separated)  $ \hbar $--adic  completion
of} \;\; \kh \otimes_A \uhat_q(\gerg^*) \; =: \; U_\hbar(\gerg^*)
\, ; \quad  \cr }  $$   
this and (2.5) together give  
  $$  {U_\hbar(\gerg)}' = \, \hbox{(separated)  $ \hbar $--adic 
completion of } \; \kh \otimes_A \utilde_q(\gerg) \, .  
\eqno (2.6)  $$  
This gives us a concrete description of  $ \, {U_\hbar(\gerg)}'
\, $: if  $ U_\hbar(\gerg) $  is topologically generated   --- as
usual ---  by Chevalley-like generators  $ \, F_i $,  $ H_j $, 
$ E_i $  (for $ i $  and  $ j $  in some set of indices  $ I $ 
and  $ J $, depending on the type of  $ \gerg $)  and if the 
$ F_\alpha $'s, resp.~$ E_\alpha $'s,  are (quantum) root
vectors attached to the positive, resp.~negative, roots
of  $ \gerg $  (like, for instance, in [Ga3--4]), then  
                                                \par   
   {\it  $ \, {U_\hbar(\gerg)}' \, $  is the unital topological
subalgebra of  $ U_\hbar(\gerg) $  topologically generated by the 
set  $ \, {\big\{ \dot{F}_\alpha, \dot{E}_\alpha \big\}}_\alpha 
\bigcup {\big\{ \dot{H}_j \big\}}_j \; $  with  $ \, 
\dot{F}_\alpha := \hbar F_\alpha $,  $ \, \dot{E}_\alpha := \hbar 
E_\alpha $,  $ \, \dot{H}_j := \hbar H_j \, $  for all  $ \alpha $  
and all  $ j \, $. }  
                                                \par   
   Having this description in our hands, we can recognize that 
Theorem 2.3 in this case is also proved in [Ga1], Theorem 4.4 (or 
simply Corollary 3.8, for  $ \, c = 1 \, $), for the finite case, 
and in [Ga2],  Corollary 2.5{\it (b)},  for the affine case.  

\vskip2,1truecm

\centerline {\bf \S\; 3. \  Braidings from geometric quantisation: 
Weinstein and Xu's approach } 

\vskip13pt 

  {\bf 3.1 The (global) classical  $ \lagR $--matrix  (cf.~[WX]).}
\, In this section we recall from [WX] the construction of the
global  $ \lagR $--matrix  and point out how it provides a
braiding.   
                                                  \par   
   From now on, let  $ \, \Bbbk \in \big\{ {\Bbb R}, {\Bbb C}
\big\} \, $.  Let  $ \, (\gerg,\r) \, $  be a (finite dimensional)
quasitriangular Lie bialgebra, and write $ \, \r = \sum_i r_i^+
\otimes r_i^- \in \gerg \otimes \gerg \, $. Define linear maps   
  $$  \r_\pm \, \colon \, \gerg^* \loongrightarrow \gerg^{**} =
\gerg \; ,  \qquad  \r_\pm(\eta):= \pm \sum_i \eta(\r_i^\pm)
\cdot r_i^\mp \qquad  \forall \; \eta \in \gerg^* \, .  
\eqno (3.1)  $$   
These are both Lie algebra homomorphisms; if  $ (\gerg,\r) $ 
is triangular, then  $ \, \r_+ = \r_- \, $.  
                                                \par  
   Let  $ G $  be a complete Poisson Lie group, and assume a dual
Poisson Lie group  $ G^* $  exists (in general,  {\sl only a 
germ\/} of such a group is defined); then their tangent Lie 
bialgebras $ \gerg $  and  $ \gerg^* $  are dual to each other.  
We say that $ G $  {\sl is quasitriangular\/}  if  $ \gerg $  is 
quasitriangular {\sl and\/}  if the Lie algebra homomorphisms  $ 
\, \r_\pm \, \colon \, \gerg^* \longrightarrow \gerg \, $  defined 
above  {\sl lift\/}  to Lie group homomorphisms  $ \, R_\pm \, 
\colon \, G^* \longrightarrow G \, $.  In this case, we define    
  $$  \phi \, ,  \, \psi \, \colon \, G^* \loongrightarrow G \; , 
\qquad  \phi(x):= R_+(x^{-1}) \, ,  \quad \psi(x):= R_-(x^{-1}) \, 
, \qquad  \forall \; x \in G^* \; .   \eqno (3.2)  $$   These are 
both Poisson morphisms; if  $ G $  is triangular (i.e.~the like is 
true for  $ \gerg $)  then  $ \, R_+ = R_- \, $,  \, hence $ \, 
\phi = \psi \, $.  
                                                \par  
   We shall use the following conventions for  {\sl dressing
transformations\/}:  the left and right dressing transformation of 
$ G $  on  $ G^* $  are denoted, respectively, by  $ \, 
\lambda_g{}u \, $  and  $ \, \rho_g{}u \, $  for all  $ \, g \in G 
\, $  and  $ \, u \in G^* \, $.  Similarly, we denote the left and 
right dressing transformation of  $ G^* $  on  $ G $  by  $ \, 
\lambda_u{}g \, $ and  $ \, \rho_u{}g \, $  for all  $ \, u \in 
G^* \, $  and  $ \, g \in G \, $.  
                                                \par  
   By definition, the  {\it global classical  $ \lagR $--matrix\/} 
is   
  $$  \lagR := \Big\{ \big(\psi\big(v^{-1}\big), u,
\phi(\lambda_{\psi(v^{-1})}u), v\big) \,\Big|\; u, v \in G^* 
\Big\} = \Big\{ \big(\psi\big(v^{-1}\big), u, \rho_{v^{-1}} 
\phi(u), v\big) \,\Big|\; u, v \in G^* \Big\}  $$   
% 
% 
%   $$  \hbox{ $ \eqalign{ 
%    \Cal{R}  &  := \Big\{\, \big(\psi\big(v^{-1}\big), u,
% \phi(\lambda_{\psi(v^{-1})}u), v\big) \;\Big|\;\, u,
% v \in G^* \Big\} =  \cr  
%    &  \phantom{:}= \Big\{\, \big(\psi\big(v^{-1}\big), u,
% \rho_{v^{-1}} \phi(u), v\big) \;\Big|\;\, u, v \in G^* \Big\} 
% \quad  \big( \subseteq D \times D \big)  \cr } $ }   \eqno (3.3)  $$   
% 
% 
which is a Lagrangian submanifold of  $ D \times D $.  It is shown 
in [WX] how this object enjoys a bunch of properties which are 
exactly the analogous of those of a quantum  $ R $--matrix;  in
addition, if  $ G $ is  {\sl triangular\/},  then  $ \lagR $  is 
{\sl unitary},  by which we mean that  $ \, \lagR^{\text{op}} =
\lagR^{-1} \, $  (in the sense of [WX], Remark 8.3).  Moreover,
these properties imply the following result:   

\vskip9pt  

\proclaim{Theorem 3.2} (cf.~[WX], Corollary 7.2) If  $ G $  is a 
complete quasitriangular Poisson Lie group, then the map  $ \; 
\Cal{R} = \Cal{R}_{{}_{W\!X}} \, \colon \, G^* \hskip-1,5pt \times 
G^* \loongrightarrow G^* \hskip-1,5pt \times G^* \; $  given by  
  $$  (u,v) \mapsto \big( \lambda_{\psi(v^{-1})}u \, , \,
\lambda_{\phi(\lambda_{\psi(v^{-1})}u)}v \big) = \big(
\lambda_{\psi(v^{-1})}u \, , \, \rho_{\phi(u^{-1})}v \big) 
\qquad  \forall \;\; u, v \in G^*  $$   
is a Poisson diffeomorphism such that   
  $$  \hbox{ $ \eqalign{ 
   m \, \circ \,  &  \Cal{R} = m^{\text{op}}  \cr  
   \Cal{R} \circ (m \otimes \hbox{\rm id}\,) = (m \otimes
\hbox{\rm id}\,) \circ \Cal{R}_{2{}3} \circ \Cal{R}_{1{}3} \; ,  & 
\qquad  \Cal{R} \circ (\hbox{\rm id} \otimes m) = (\hbox{\rm id}
\otimes m) \circ \Cal{R}_{1{}2} \circ \Cal{R}_{1{}3}  \cr } $ }  
\eqno(3.3)  $$   
where  $ m $  is the product of  $ G^* $  and  $ \, m^{\text{op}}
:= m \circ \sigma \, $  (with  $ \sigma $  as in \S 1.8).  In
particular, $ \Cal{R} $  is a solution of the QYBE, and it
restricts to a similar mapping  $ \, \Cal{S} \times \Cal{S}
\loongrightarrow \Cal{S} \times \Cal{S} \, $  for every symplectic
leaf  $ \, \Cal{S} $   of  $ G^* $.  In addition, if  $ G $  is
triangular then  $ \Cal{R} $  is  {\sl unitary}, which means 
$ \, \Cal{R}^{-1} = \sigma \circ \Cal{R} \circ \sigma \, $.  
\endproclaim  

\demo{Proof}  It is just a matter of recalling or reformulating
some results of [WX].  The identity in the first line of (3.3)
is proved by Theorem 5.1 in  [{\it loc.~cit.}];  the second line
of identities instead is a simple reformulation of Theorem 5.4 in 
[{\it loc.~cit.}];  finally, in the triangular case the unitarity
of  $ \Cal{R} $  follows from the unitarity of the global 
$ \lagR $--matrix  $ \lagR $,  by Corollary 8.2 and  {\it
Remark 8.3\/}  in  [{\it loc.~cit.}].   \qed  
\enddemo  

\vskip9pt  

\proclaim{Corollary 3.3}  The mapping   
  $$  \R_{{}_{W\!X}} := \Cal{R}^* \, \colon \, F[G^*] \otimes
F[G^*] = F[G^* \hskip-1,5pt \times G^*] \llongrightarrow
F[G^* \hskip-1,5pt \times G^*] = F[G^*] \otimes F[G^*]  $$   
naturally induced by  $ \Cal{R} $  is a braiding, which is unitary
if  $ \, \Cal{R} $  is.  In particular, this canonically induces a
braiding  $ \; \R_{{}_{W\!X}} \, \colon F [[\gerg^* \hskip-1,5pt
\oplus \gerg^*]] \llongrightarrow F[[\gerg^* \hskip-1,5pt \oplus
\gerg^*]] \, . $   Furthermore, the asso\-ciated infinitesimal
braiding  $ \; \overline\R_{{}_{W\!X}} \, \colon \gerg \oplus \gerg
\loongrightarrow \gerg \oplus \gerg \; $ (cf.~\S 2.4) is trivial,
i.e.~$ \, \overline\R_{{}_{W\!X}} = \hbox{\rm id}_{\gerg \oplus
\gerg} \, $.  
\endproclaim 

\demo{Proof}  The first part of the claim   --- 
$ \R_{{}_{W\!X}} $  being a braiding, unitary if 
$ \Cal{R}_{{}_{W\!X}} $  is ---   follows trivially
from Theorem 3.2 by duality; then  $ \R_{{}_{W\!X}} $ 
automatically induces an infinitesimal braiding 
$ \overline\R_{{}_{W\!X}} $  as well.   
                                              \par    
   To prove the second part   --- that is, 
$ \overline\R_{{}_{W\!X}} $  being trivial ---  
we must go back to the definition and the properties
of dressing actions.  Recall that the left dressing action
of  $ G $  on  $ G^* $  is defined as follows.  For all  $ \, g
\in G $,  $ \gamma \in G^* $,  there exist unique  $ \, g^\gamma
\in G $,  $ \gamma^g \in G^* \, $  such that  $ \; g \cdot \gamma
= \gamma^g \cdot g^\gamma \; $;  \; then the left dressing action 
$ \; \lambda \, \colon \, G \times G^* \! \loongrightarrow G^*
\; $  of  $ G $  on  $ G^* $  is given by  $ \; \lambda_g(\gamma)
\equiv \lambda(g,\gamma) := \gamma^g \, $,  \; for all  $ \, g
\in G $,  $ \gamma \in G^* $.  
                                              \par    
   Now, for all  $ \, X \in \gerg $,  $ Y \in \gerg^* $  and 
$ t \in {\Bbb R} $,  \, we have   
  $$  \exp(t \, X) \cdot \exp(t \, Y) = {\exp(t \, Y)}^{\exp(t
\, X)} \cdot {\exp(t \, X)}^{\exp(t \, Y)}  $$   
whence Taylor expansion gives   
  $$  \displaylines{ 
   \qquad   {\exp(t \, Y)}^{\exp(t \, X)} = {\big( 1 + t \, Y
+ t^2 \, Y^2 \big/ 2 + \cdots \big)}^{( 1 \, + \, t X \, + \,
t^2 X^2 / 2 \, + \, \cdots )} =   \hfill  \cr  
   \hfill   = 1 + t \, Y + t^2 \, Y^X + \cdots + t^2 \,
Y^2 \big/ 2 + \cdots  \qquad  \cr }  $$   
(where  $ \, Y^X \, $  denotes the action of  $ X $  onto  $ Y $ 
induced at the infinitesimal level by the dressing action), hence
at first order in  $ t \, $  we have simply  $ Y \, $!  Applied to
the situation  $ \, \exp(t \, X) = \psi(v^{-1}) $,  $ \exp(t \, Y)
= u \, $  this says that the first entry of  $ \, T_{(e,e)}
(\Cal{R}_{{}_{W\!X}})(Y,V) $  is just  $ Y \, $  (here  $ \,
V := \log(v) \, $,  \, and  $ e \, $  denotes the unit element
of  $ G^* \, $).  Similarly, carrying out a like analysis on
the right dressing action we get that the second entry of 
$ \, T_{(e,e)}(\Cal{R}_{{}_{W\!X}})(Y,V) $  is simply  $ V $. 
Therefore,  $ \, T_{(e,e)}(\Cal{R}_{{}_{W\!X}}) = \hbox{\rm
id}_{\gerg^* \oplus \gerg^*} \, $;  \; as  $ \,
\overline\R_{{}_{W\!X}} \, $  is just the dual of 
$ T_{(e,e)}(\Cal{R}_{{}_{W\!X}}) $,  it is trivial
as well, q.e.d.   \qed  
\enddemo  

\vskip9pt 

  {\bf 3.4 The factorizable case.} \, Let  $ \, (\gerg,\r) \, $  
be a quasitriangular Lie bialgebra: if the bilinear form on  $ \, 
\gerg \otimes \gerg \, $  naturally associated to  $ \, \r + 
\r^{\text{op}} \, $  is non-degenerate, then  $ \, 
(\gerg,\r) \, $  is said to be  {\sl factorizable}.  In this case, 
the corresponding linear map  $ \; j:= \r_+ - \r_- \, \colon \, 
\gerg^* \!\loongrightarrow \gerg \; $  is invertible. Now let 
$ G $  be a Poisson Lie group corresponding to the bialgebra above,
and let  $ G^* $  be its connected, simply connected Poisson dual.  
The Lie algebra morphisms  $ \; \r_\pm \, \colon \, \gerg^* \!
\loongrightarrow \gerg \; $  lift to group morphisms  $ \; R_\pm
\, \colon \, G^* \!\loongrightarrow G \, $, \; thus we may define 
the map  $ \; J \, \colon \, G^* \! \loongrightarrow G \; $  by 
$ \; J(u) := R_+(u) \, {R_-(u)}^{-1} \; $  (for all  $ \, u \in
G^* \, $)  whose derivative at the identity element  $ \, u \in
G^* \, $  is  $ j \, $  (note that neither  $ j $  nor  $ J $  is
a morphism).  When  $ J $  is a global diffeomorphism, we say
that  {\sl the group  $ G $  is factorizable},  since for each 
$ \, g \in G \, $  we have the factorization  $ \; g = g_+ \, 
{g_-}^{\hskip-5pt -1} \, $,  \; where  $ \; g_\pm := R_\pm \big( 
J^{-1}(g) \big) \, $.  Thanks to [WX], Proposition 9.1, any
connected, simply connected, factorizable Poisson Lie group
is complete.  
                                                 \par   
   Now, factorizability enables us to describe the classical 
$ \Cal{R} $--matrix  quite explicitly:  

\vskip9pt  

\proclaim{Theorem 3.5} (cf.~[WX], Theorem 9.2) \, Let  $ G $  be a 
factorizable Poisson Lie group, and use  $ \; J \, \colon \, G^* 
\! \loongrightarrow G \; $  to identify  $ G^* $  with  $ G $ 
(hence also  $ \, G \hskip-0,5pt \times G \, $  with  $ \, G^* 
\hskip-1,5pt \times G^* \, $).  Then:   
                                      \hfill\break   
   \indent   (a) \; the (global) classical  $ \lagR $--matrix  $ \,
\lagR \in (G \times G) \times (G \times G) \, $  takes the form   
  $$  \lagR = \Big\{\, \Big( y_- \, , \, x \, , {{\big( y_- \,
x \, {y_-}^{\hskip-5pt -1} \big)}_{\!+}}^{\hskip-5pt -1}, \,
y \Big) \;\Big|\; \forall \; x, y \in G \,\Big\} \; ;  $$  
   \indent   (b) \; the map  $ \; \Cal{R} = \Cal{R}_{{}_{W\!X}}
\, \colon \, G \hskip-1,5pt \times G \loongrightarrow G
\hskip-1,5pt \times G \; $  of Theorem 3.2 above is given by   
   $$  \Cal{R}_{{}_{W\!X}}(x,y) = \Big( y_- \, x \,
{y_-}^{\hskip-5pt -1}, {{\big( y_- \, x \, {y_-}^{\hskip-5pt -1} 
\big)}_+}^{\hskip-5pt -1} \, y \, {{\big( y_- \, x \, 
{y_-}^{\hskip-5pt -1} \big)}_+}^{\hskip-5pt +1} \Big)   \qquad  
\forall \;\; (x,y) \in G \times G \; .   \qed  $$   
\endproclaim  

\vskip9pt  

   {\bf Remark 3.6:} \, As we pointed out in the Introduction,
one can carry over the construction of Weinstein and Xu in
purely local terms, just performing it on the germ of
Poisson group underlying the quasitriangular Lie bialgebra 
$ (\gerg,\r) $,  and eventually get a braiding  $ \;
\R_{{}_{W\!X}} \, \colon F [[\gerg^* \hskip-1,5pt \oplus
\gerg^*]] \llongrightarrow F[[\gerg^* \hskip-1,5pt \oplus
\gerg^*]] \; $  and an associated infinitesimal braiding 
$ \; \overline\R_{{}_{W\!X}} \, \colon \gerg \oplus \gerg
\loongrightarrow \gerg \oplus \gerg \, $.  Our next result
is that the latter is always trivial whenever  $ (\gerg,\r) $ 
is factorizable.   

\vskip9pt  

\proclaim{Proposition 3.7} \, Let the quasitriangular Lie
bialgebra  $ (\gerg,\r) $  be factorizable.  Then the
infinitesimal braiding  $ \; \overline\R_{{}_{W\!X}} \,
\colon \gerg \oplus \gerg \loongrightarrow \gerg \oplus
\gerg \; $  is trivial,  i.e.~$ \, \overline\R_{{}_{W\!X}}
= \hbox{\rm id}_{\gerg \oplus \gerg} \, $.   
\endproclaim  

\demo{Proof}  Let  $ G_{\text{\it loc}} $  be the germ of
Poisson group associated to the Lie bialgebra  $ \gerg $. 
Then the "local" version of  Theorem 3.5{\it (b)}  ensures
that the map  $ \; \Cal{R}_{{}_{W\!X}} \, \colon \, G_{\text{\it
loc}} \hskip-1,5pt \times G_{\text{\it loc}} \! \loongrightarrow
G_{\text{\it loc}} \hskip-1,5pt \times G_{\text{\it loc}} \; $ 
is given by  $ \; \Cal{R}_{{}_{W\!X}}(x,y) = \Big( y_- \, x \,
{y_-}^{\hskip-5pt -1}, {{\big( y_- \, x \, {y_-}^{\hskip-5pt
-1} \big)}_+}^{\hskip-5pt -1} \, y \, {{\big( y_- \, x \,
{y_-}^{\hskip-5pt -1} \big)}_+}^{\hskip-5pt +1} \Big) \; $ 
for all  $ \; x, y \in G_{\text{\it loc}} \, $.  Now, for
all  $ \, A $,  $ B \in \gerg \, $  and  $ \, t \in {\Bbb R}
\, $  we have   
  $$  \displaylines{ 
   {} \quad  \exp(t\,A) \, \exp(t\,B) \, {\exp(t\,A)}^{-1} =  
\hfill  \cr  
   = \big( 1 + t \, A + t^2 A^2 \big/ 2 + \cdots \big)
\, \big( 1 + t \, B + t^2 B^2 \big/ 2 + \cdots \big) \,
\big( 1 - t \, A + t^2 A^2 \big/ 2 - \cdots \big) =  \cr   
   {} \hfill   = \, 1 \, + \, t \, B \, + \, t^2 \big(
2 \, (AB-BA) \, + \, B^2 \big) \big/ 2 \, + \, \cdots \; ; 
\quad  \cr }  $$ 
applying this recipe to  $ \, A = \log(y_-) $,  $ \, B =
\log(x) $,  \, and looking at first order (in  $ t \, $) 
we find out that the first entry of  $ \, T_{(e,e)}
(\Cal{R}_{{}_{W\!X}})(x,y) $  is just  $ x \, $;  \,
similarly we get that the second entry of  $ \, T_{(e,e)}
(\Cal{R}_{{}_{W\!X}})(x,y) $  is  $ y $.  Thus  $ \, T_{(e,e)}
(\Cal{R}_{{}_{W\!X}}) \, $  is the identity, and since  $ \,
\overline\R_{{}_{W\!X}} \, $  is just its dual, it is the
identity as well, q.e.d.   \qed   
\enddemo   

\vskip2,1truecm

\centerline {\bf \S\; 4. \  Comparing the braidings  $ \;
\R_{{}_{W\!X}} \, $  and  $ \; \R_{{}_{G\!H}} \, $:  \,
the case of  $ \, \gerg = \gersl_2 \, $. } 

\vskip13pt 

{\bf 4.1 The general problem.} \, We noticed that the construction
of [WX] can be performed for any quasitriangular Lie bialgebra by
acting locally, so to get a braiding  $ \; \R_{{}_{W\!X}} \, $  on
the dual formal Poisson group, exactly like one can do following
[GH] to get a braiding  $ \; \R_{{}_{G\!H}} \, $.  Since these
braidings share similar properties   --- like functoriality and
infinitesimal triviality, for instance ---   we are led to raise
the following   

\vskip5pt   

\proclaim{Question}  Given any quasitriangular Lie bialgebra 
$ \gerg $,  do the braidings  $ \; \R_{{}_{W\!X}} \, $  and  $ \;
\R_{{}_{G\!H}} \, $  on  $ \, F[[\gerg^*]] \, $  coincide?   
\endproclaim   

\vskip5pt   

   The purpose of the present section is to provide a positive
answer to this question for the simplest case of  $ \, \gerg
= \gersl_2({\Bbb C}) \, $.   

\vskip9pt   

  {\bf 4.2 The geometrical setting.} \, In this section, let  $ \,
\Bbbk = {\Bbb C} \, $.  Let  $ \, G := {SL}_2 \equiv {SL}_2
({\Bbb C}) \, $. Its tangent Lie algebra  $ \, \gerg = \gersl_2
\, $  is generated by  $ \, f $,  $ h $,  $ e \, $  (the  {\it
Chevalley generators\/})  with relations  $ \, [h,e] = 2 \, e $, 
$ [h,f] = -2 \, f $,  $ [e,f] = h \, $.  The formul\ae{}  $ \,
\delta(f) = (f \otimes h - h \otimes f) \big/ 2 \, $,  $ \,
\delta(h) = 0 \, $,  $ \, \delta(e) = (e \otimes h - h \otimes e)
\big/ 2 \, $,  define a Lie cobracket on  $ \gerg \, $:  \, indeed,
this makes  $ \gersl_2 $  into a quasitriangular Lie bialgebra,
whose  $ \r $-matrix  is  $ \; \r := e \otimes f + (h \otimes h)
\big/ 4 \, $.  This corresponds to a structure of complex Poisson
Lie (actually,  {\sl algebraic\/}) group on  $ G $,  which is
complete and quasitriangular.  
                                             \par
% 
% 
%   The function algebra  $ F[{SL}_2] $  is the unital associative
% commutative  $ \Bbbk $--algebra  with generators  $ \, a $,  $ b $, 
% $ c $,  $ d \, $  and the relation  $ \, a d - b c = 1 \, $,  and
% Poisson Hopf structure given by  
%   $$  \displaylines{
%   \Delta(a) = a \otimes a + b \otimes c \, ,  \;\,
% \Delta(b) = a \otimes b + b \otimes d \, ,  \;\,
% \Delta(c) = c \otimes a + d \otimes c \, ,  \;\,
% \Delta(d) = c \otimes b + d \otimes d  \cr
%   \epsilon(a) = 1  \, ,  \hskip5pt  \epsilon(b) = 0 \, ,  \hskip5pt
% \epsilon(c) = 0 \, ,  \hskip5pt  \epsilon(d) = 1 \, ,  \hskip13pt 
% S(a) = d \, ,  \hskip5pt  S(b) = -b \, ,  \hskip5pt  S(c) = - c \, ,
% \hskip5pt  S(d) = a  \cr
%   \{a,b\} = b a \, , \quad  \{a,c\} = c a \, ,  \quad  \{b,c\}
% = 0 \, ,  \quad  \{d,b\} = - b d \, ,  \quad  \{d,c\} = - c d \, ,
% \quad  \{a,d\} = 2 b c \, .  \cr }  $$  
%                                              \par
% 
% 
   In the dual Lie bialgebra  $ \, \gerg^* = {\gersl_2}^{\!*} \, $, 
\, let  $ \, \big\{ e^*, f^*, h^* \big\} \, $  be the basis dual 
to $ \, \{e,f,h\} \, $,  \, and consider the basis  $ \, \big\{ \e 
:= e^*, \f := f^*, \h := - 2 \, h^* \big\} \, $.  Then the Lie 
bialgebra structure of  $ \, {\gersl_2}^{\!*} \, $  is described 
by the formul\ae{}  $ \, [\h,\e] = \e $,  $ [\h,\f\,] = \f $,  $ 
[\e,\f\,] = 0 \, $,  and  $ \, \delta(\f\,) = \h \otimes \f - \f 
\otimes \h $, $ \, \delta(\h) = 2 \, (\f \otimes \e - \e \otimes
\f \,) $,  $ \, \delta(\e) = \e \otimes \h - \h \otimes \e \, $. 
Then  $ {\gersl_2}^{\!*} \, $  can be realized as the Lie
algebra of pairs of matrices   
  $$  {\gersl_2}^{\!*} = \Bigg\{\,
\bigg( \bigg( \matrix  -t & 0 \\  c & t  \endmatrix \bigg) , 
\bigg( \matrix  t & b \\  0 & -t  \endmatrix \bigg) \bigg) 
\,\Bigg\vert\, b, c, t \in \Bbbk \,\Bigg\} \; \subseteq \; \gersl_2 
\hskip-0,5pt \times \gersl_2   \eqno (4.1)  $$   
(with the Lie subalgebra structure inside  $ \, \gersl_2
\hskip-0,5pt \times \gersl_2 \, $).  It follows that the unique
connected simply connected complex Poisson Lie (actually, 
{\sl algebraic\/})  group whose tangent Lie bialgebra is 
$ {\gersl_2}^{\!*} \, $  can be realized as the group of
pairs of matrices (the left subscript  $ s $  meaning
"simply connected")   
  $$  {}_s{{SL}_2}^{\!*} = \Bigg\{\,
\bigg( \bigg( \matrix  z^{-1} & 0 \\  y & z  \endmatrix \bigg) , 
\bigg( \matrix  z & x \\  0 & z^{-1}  \endmatrix \bigg) \bigg) 
\,\Bigg\vert\, x, y \in k, z \in \Bbbk \setminus \{0\} \,\Bigg\} \; 
\subseteq \; {SL}_2 \hskip-0,5pt \times {SL}_2   \eqno (4.2)  $$ 
(with the subgroup structure inside  $ \, {SL}_2 \hskip-0,5pt 
\times {SL}_2 \, $);  this group has a "small" centre, namely  $ 
\, Z := \big\{ (I,I\,), (-I,-I\,) \big\} \, $,  so there is only 
one other (Poisson) group sharing the same Lie (bi)algebra, namely 
the quotient $ \, {}_a{{SL}_2}^{\!*} := {}_s{SL_2}^* \Big/ Z \, $  
(the adjoint of $ \, {}_s{{SL}_2}^{\!*} \, $,  as the left 
subscript  $ a $  means). Therefore  $ F\big[ {}_s{{SL}_2}^{\!*} 
\big] $  is the unital associative commutative  $ \Bbbk $--algebra  
with generators  $ \, x $, $ z^{\pm 1} $,  $ y $,  with Poisson 
Hopf structure given by   
  $$  \displaylines{ 
   \Delta(x) = x \otimes z^{-1} + z \otimes x \, ,  \hskip21pt 
\Delta\big(z^{\pm 1}\big) = z^{\pm 1} \otimes z^{\pm 1} \, , 
\hskip21pt  \Delta(y) = y \otimes z^{-1} + z \otimes y  \cr 
  \epsilon(x) = 0  \, ,  \hskip10pt  \epsilon\big(z^{\pm 1}\big) = 1
\, ,  \hskip10pt  \epsilon(y) = 0 \, ,  \hskip31pt 
  S(x) = -x \, ,  \hskip10pt  S\big(z^{\pm 1}\big) = z^{\mp 1} \, , 
\hskip10pt  S(y) = -y  \cr 
  \{x,y\} = z^{-2} - z^{+2} \, ,  \hskip27pt 
\big\{z^{\pm 1},x\big\} = \mp x z^{\pm 1} \big/ 2 \, ,  \hskip27pt 
\big\{z^{\pm 1},y\big\} = \pm z^{\pm 1} y \big/ 2  \cr }  $$  
(N.B.: with respect to this presentation, we have  $ \, \text{f} = 
{\partial_y}{\big\vert}_u \, $,  $ \, \text{h} = {\, z \, \over \, 
2 \,} \, {\partial_z}{\big\vert}_u \, $,  $ \, \text{e} = 
{\partial_x}{\big\vert}_u \, $,  where  $ u $  is the identity 
element of  $ {}_s{SL_2}^* \, $).  Moreover,  $ F \big[ 
{}_a{{SL}_2}^{\!*} \big] $ can be identified with the Poisson Hopf 
subalgebra of  $ F\big[{}_s{SL_2}^*\big] $  spanned by products of 
an even number of generators, i.e.~monomials of even degree: this 
is generated as a unital subalgebra, by $ x z $,  $ z^{\pm 2} $, 
and  $ z^{-1} y $.  Finally, the (algebra of regular functions on
the) Poisson algebraic formal group  $ F\big[\big[ {\gersl_2}^{\!*}
\big]\big] $  is the  $ \Ker\,(\epsilon) $--adic  completion
of both  $ F \big[{}_s{{SL}_2}^{\!*}\big] $  and  $ F \big[
{}_a{{SL}_2}^{\!*}\big] $;  in the first case  $ \Ker\,
(\epsilon) $  is generated (as an ideal) by  $ \, x \, $, 
$ \, \big( z^{\pm 1} - 1 \big) \, $  and $ \, y \, $,  \,
therefore  $ \; F \big[\big[ {\gersl_2}^{\!*} \big]\big] =
\Bbbk \big[\big[ x, (z-1), y \big]\big] \; $  as a topological 
$ \Bbbk $--algebra  (note that  $ \, z^{-1} - 1 = \sum_{n>0}
{(-1)}^n {(z-1)}^n \, $,  \, so the generator  $ \, z^{-1}
- 1 \, $  is superfluous) with the unique Poisson Hopf
structure which extends by continuity the one on 
$ F \big[{}_s{{SL}_2}^{\!*}\big] $.   

\vskip9pt 

  {\bf 4.3 Weinstein and Xu's construction.} \, In the
framework of \S 4.2, let  $ \, G := {SL}_2 \, $,  $ \, G^*
:= {}_s{{SL}_2}^{\!*} \, $.  In this section we compute the
braiding  $ \R_{{}_{W\!X}} $  for  $ G \, $;  \, despite the
fact that not all requirements of [WX] are fulfilled, we can
show that that construction can still be carried out at the
local level: to fulfill this goal is then just a matter of
matrix computation.   
                                              \par     
   It follows from definitions   --- cf.~[WX], \S 9 ---   that
the maps  $ \, \r_\pm \, \colon \, \gerg^* \loongrightarrow \gerg
\, $  are given by  $ \; \r_+(\e) = 2 \, f \, $,  $ \, \r_+(\h) =
-h \big/ 2 \, $,  $ \, \r_+(\f\,) = 0 \, $,  $ \, \r_-(\e) = 0
\, $, $ \, \r_-(\h) = +h \big/ 2 \, $,  $ \, \r_-(\f) = -2 \,
\e \, $,  \; and the maps  $ \; R_\pm \, \colon \, G^*
\hskip-1,5pt \longrightarrow G \; $  are, respectively,
the projection to the second and the first factor w.r.t.~to
the description of  $ \, G^* = {}_s{{SL}_2}^{\!*} \, $  in
(4.2).  Then for the maps  $ \; j \, \colon \, \gerg^* \! 
\loongrightarrow \gerg \; $  and  $ \; J \, \colon \, G^* \! 
\loongrightarrow G \; $  defined in \S 3.4 we have that  $ j $
is bijective but  $ J $  is not, for it has kernel  $ \, \Ker
\,(J\,) = Z \, $  (hence it is a 2--to--1 map) and image    
  $$  \hbox{\rm Im}\,(J\,) = G^0 := \left\{\, 
   \pmatrix  a  &  b  \\   c  &  d   \endpmatrix  \;\bigg|\;\;
a, b, c, d \in \Bbb{C}, \; d \not= 0 \;\right\}  $$   
that is the  {\sl big cell\/}  of  $ \, G = {SL}_2 \, $:  \, in
fact,  $ J $  is an unramified 2-fold covering of  $ G^0 $. 
Therefore,  $ J $  is not a global diffeomorphism, but it
factors to a global diffeomorphism   
  $$  J_a \, \colon \, {}_aG^* \equiv
{}_a{{SL}_2}^{\!*} := {}_s{SL_2}^* \Big/ Z  \;\;\;{\buildrel 
\simeq \over 
{\joinrel\relbar\joinrel\relbar\joinrel\llongrightarrow}}\;\;
G^0  $$   
given by  $ \; J_a\big( g \cdot Z \big) := J(g) \; $  for all 
$ \, g \in {}_s{SL_2}^* \, $.  We need a section of  $ J $  and
of  $ J_a \, $.  Since   
  $$  \displaylines{ 
   J \left( \hskip-2pt  \pmatrix  A^{-1}  &  0  \\  
                                       B  &  A^{+1}  \endpmatrix , 
\pmatrix  A^{+1}  &  C  \\   0  &  A^{-1}  \endpmatrix 
\hskip-2pt \right)  =   \hfill  \cr  
   \hfill   = \pmatrix  A^{+1}  &  C  \\   0  &  A^{-1} 
\endpmatrix  \cdot  \pmatrix  A^{+1}  &  0  \\   -B  &  A^{-1} 
\endpmatrix  =  \pmatrix  A^{+2} - B \, C  &  A^{-1} C  \\ 
- A^{-1} B  &  A^{-2}  \endpmatrix  \cr }  $$   
we have  $ \; J \left( \hskip-2pt  \pmatrix  A^{-1}  &  0  \\ 
B  &  A^{+1}  \endpmatrix , \pmatrix  A^{+1}  &  C  \\   0  & 
A^{-1}  \endpmatrix  \hskip-2pt \right)  =  \pmatrix  a  & 
b  \\   c  &  d  \endpmatrix \; $  if and only if   
  $$  A = \pm d^{-1/2} \, ,  \qquad  B = \pm b \, d^{-1/2} \, ,  
\qquad  C = \mp c \, d^{-1/2}   \eqno (4.3)  $$   
for any matrix  $ \; \pmatrix  a  &  b  \\   c  &  d  \endpmatrix
\in G^0 \, $;  \, these formul\ae{} clearly define two
differentiable sections of  $ J $  (taking either upper
or lower signs) and one of $ J_a $  (for which the sign
is irrelevant).  

\vskip7pt 

   {\it Remark:} \, although  $ G $ is not factorizable,
nevertheless we can still use  Theorem 3.5{\it (b)\/} 
to compute the map  $ \Cal{R}_{{}_{W\!X}} $,  namely   
   $$  \Cal{R}_{{}_{W\!X}}\big(X',Y'\big) = \Big( J^{-1} \Big( Y_-
\, X \, {Y_-}^{\hskip-5pt -1} \Big), J^{-1} \Big( {{\big( Y_- \, X 
\, {Y_-}^{\hskip-5pt -1} \big)}_+}^{\hskip-5pt -1} \, Y \, {{\big( 
Y_- \, X \, {Y_-}^{\hskip-5pt -1} \big)}_+}^{\hskip-5pt +1} \Big) 
\Big)   \eqno (4.4)  $$   
for all  $ \; \big(X',Y'\big) \in G^* \hskip-1,5pt \times G^* \, $ 
and  $ \, (X,Y) := \big( J(X'), J(Y') \big) \in G \times G \, $, 
\, where  $ J^{-1} $  is one of the two aforesaid sections of 
$ J $,  namely the unique one such that the resulting 
$ \Cal{R}_{{}_{W\!X}}(X',Y') $  map  $ (e_{{}_{G^*}},
e_{{}_{G^*}}) $ onto itself.  In fact, although  $ J $ is not a 
diffeomorphism it is nevertheless a (finite) covering on  $ G^0 $,  
hence it is a {\sl local diffeomorphism\/} (around the identity 
element  $ \, e_{G^*} \in G^* \, $)  on $ G^0 $,  therefore the 
description of  $ \Cal{R}_{{}_{W\!X}} (X',Y') $  afforded by  
Theorem 3.5{\it (b)\/}, through  $ J $ and a local section  $ 
J^{-1} $,  is still available (locally around  $ \, 
(e_{{}_{G^*}},e_{{}_{G^*}}) \in G^* \hskip-1,5pt \times G^* \, $):  
\, to have a  {\sl global\/}  description, one has just to choose 
the unique section  $ J^{-1} $  which maps $ \, e_{G^0} = e_G \, $  
onto  $ e_{G^*} $.  Therefore, we shall now go on computing 
$ \Cal{R}_{{}_{W\!X}} $  following this strategy.   

\vskip7pt 

   Let  $ \; X := J \bigg( \hskip-3pt \pmatrix  z^{-1} & 0 \\ 
y & z  \endpmatrix , \pmatrix  z & x \\  0 & z^{-1}  \endpmatrix 
\hskip-3pt \bigg) \, $,  $ \, Y := J \bigg( \hskip-3pt \pmatrix 
\zeta^{-1} & 0 \\  \eta & \zeta  \endpmatrix , \pmatrix  \zeta & 
\chi \\  0 & \zeta^{-1}  \endpmatrix \hskip-3pt \bigg) \in J(G^*) 
= G^0 \, $.  Then we have   
  $$  \displaylines{ 
   Y_- \cdot X \cdot {Y_-}^{\hskip-3pt -1} = 
\bigg( \matrix  \zeta^{-1} & 0 \\  \eta & \zeta  \endmatrix \bigg) 
\,  \bigg( \matrix  z & x \\  0 & z^{-1} \endmatrix \bigg)  \, 
\bigg( \matrix  z & 0 \\  -y & z^{-1}  \endmatrix \bigg) \, \bigg( 
\matrix  \zeta & 0 \\  -\eta & \zeta^{-1}  \endmatrix \bigg) =   
\hfill  \cr  
   = \bigg( \matrix  \zeta^{-1} & 0 \\  \eta & \zeta \endmatrix \bigg)
\,  \bigg( \matrix  z^2 - x \, y  &  z^{-1} x \\  - z^{-1} y  & 
z^{-2}  \endmatrix \bigg)  \,  \bigg( \matrix  \zeta & 0 \\ -\eta 
& \zeta^{-1}  \endmatrix \bigg) =  \cr  
   \hfill  = \pmatrix  z^2 - x \, y - \eta \, \zeta^{-1} x z^{-1} 
&  \zeta^{-2} z^{-1} x \;  \\  \; \eta \, \zeta \, z^2 - \big( \eta
\, \zeta^{-2} \zeta + y \, z^{-1} \zeta^{+2} \big) \Theta^2  & 
z^{-2} \Theta^2  \endpmatrix  \cr }  $$   
with  $ \; \Theta := {\big( 1 + \eta \, x \, z \,
\zeta^{-1} \big)}^{1/2} \, $.  Using (4.3) we get from this  
  $$  \eqalign{ 
   {\big( Y_- \, X \, {Y_-}^{\hskip-5pt -1} \big)}_+  &  \, = 
\pm \pmatrix z^{+1} \Theta^{-1}  &  x \, \zeta^{-2} \Theta^{-1}  
\\  0  &  z^{-1} \Theta^{+1}  \endpmatrix  \, ,   \cr  
   {\big( Y_- \, X \, {Y_-}^{\hskip-5pt -1} \big)}_-  &  \, = 
\pm \pmatrix  z^{-1} \Theta^{+1}  &  0  \\  y \, \zeta^{+2} 
\Theta^{+1} + \eta \, \zeta \big( z^{-1} \Theta^{+1} - z^{+3} 
\Theta^{-1} \big) \;  &  \; z^{+1} \, \Theta^{-1} 
\endpmatrix  \cr }   \eqno (4.5)  $$   
which gives   
% 
%   $$  J^{-1} \big( Y_- \, X \, {Y_-}^{\hskip-5pt -1} \big) =
% \pm \left( \pmatrix  z^{-1} \Theta  &  0  \\  
% y \, \zeta^2 \Theta + \eta \, \zeta \big( z^{-1} \Theta
% - z^3 \Theta^{-1} \big)  &  z \, \Theta^{-1} 
% \endpmatrix , \pmatrix z \Theta^{-1}  &  x \, \zeta^{-2}
% \Theta^{-1}  \\  0  &  z^{-1} \Theta^{+1}  \endpmatrix
% \right)  $$   
% 
% 
  $$  J^{-1} \big( Y_- \, X \, {Y_-}^{\hskip-5pt -1} \big) =
\pm \hskip-3pt \left( \hskip-3pt \pmatrix  z^{-1} \Theta  &  0  \\  
y \, \zeta^2 \Theta + \eta \, \zeta \big( z^{-1} \Theta - z^3 
\Theta^{-1} \big) \!  &  \! z \, \Theta^{-1} \! \endpmatrix 
\hskip-2pt , \hskip-2pt  \pmatrix \! z \Theta^{-1} \!  & \! x \, 
\zeta^{-2} \Theta^{-1}  \\  \! 0  &  z^{-1} \Theta^{+1} 
\endpmatrix \hskip-3pt \right)   \hfill \hskip7pt  (4.6)  $$   
as possible preimages of  $ \, Y_- \, X \, {Y_-}^{\hskip-5pt -1} 
\, $  in  $ \, G^* \hskip-1,5pt \times G^* \, $.  This takes care 
of the first entry in the right-hand-side of (4.4).   
                                              \par   
   As for the second entry, we have (noting that the ambiguity of
sign in (4.5) is irrelevant)  
  $$  \displaylines{ 
   {{\big( Y_- \, X \, {Y_-}^{\hskip-5pt -1}
\big)}_+}^{\hskip-5pt -1} \cdot Y \cdot {{\big( Y_- \, X \, 
{Y_-}^{\hskip-5pt -1} \big)}_+}^{\hskip-5pt +1} =   \hfill  \cr  
   \hfill   =  \pmatrix  z^{-1} \Theta^{+1}  &  - x \, \zeta^{-2}
\Theta^{-1}  \\  0  &  z^{+1} \Theta^{-1}  \endpmatrix  \cdot 
\pmatrix  \zeta^{+1}  &  \chi  \\  0  &  \zeta^{-1}  \endpmatrix 
\cdot  \pmatrix  \zeta^{+1}  &  0  \\  -\eta  &  \zeta^{-1} 
\endpmatrix  \cdot  \pmatrix  z^{-1} \Theta^{+1}  &  - x \,
\zeta^{-2} \Theta^{-1}  \\  0  &  z^{+1} \Theta^{-1}  \endpmatrix 
=  \cr  
   =  \pmatrix  z^{-1} \Theta^{+1}  &  - x \, \zeta^{-2} \Theta^{-1} 
\\  0  &  z^{+1} \Theta^{-1}  \endpmatrix  \cdot 
\pmatrix  \zeta^{+2} - \eta \, \chi  &  \chi \, \zeta^{-1}  \\ - 
\eta \, \zeta^{-1}  &  \zeta^{-2}  \endpmatrix  \cdot \pmatrix  
z^{-1} \Theta^{+1}  &  - x \, \zeta^{-2} \Theta^{-1}  \\ 0  &  
z^{+1} \Theta^{-1}  \endpmatrix  =   \hfill  \cr  
   \hfill  =  \pmatrix  \zeta^{+2} - \eta \, \chi + \eta \, x \,
z^{+1} \zeta^{-3} \, \Theta^{-2}  &  \quad  x \, z^{-1} + \chi \, 
z^{-2} \zeta^{-1} - x \, z^{-1} \zeta^{-4} \Theta^{-2}  \\  - \eta 
\, z^{+2} \zeta^{-1} \Theta^{-2}  &  \zeta^{-2} \Theta^{-2}  
\endpmatrix \; .  \cr }  $$   
Again using (4.3) we find   
  $$  \eqalign{ 
   \Big( {{\big( Y_- \, X \, {Y_-}^{\hskip-5pt -1}
\big)}_+}^{\hskip-5pt -1}  &  \cdot Y \cdot {{\big( Y_- \, X \, 
{Y_-}^{\hskip-5pt -1} \big)}_+}^{\hskip-5pt +1} {\Big)}_+  =  \cr  
   =  \; \pm  &  \pmatrix  \zeta^{+1} \Theta^{+1}  \quad  & 
\chi \, z^{-2} \Theta^{+1} + x \, z^{-1} \zeta^{+1} \Theta^{+1} - 
x \, z^{-1} \zeta^{-3} \Theta^{-1}  \\  0  &  \zeta^{-1} 
\Theta^{-1} 
\endpmatrix   \cr  
   \Big( {{\big( Y_- \, X \, {Y_-}^{\hskip-5pt -1}
\big)}_+}^{\hskip-5pt -1}  &  \cdot Y \cdot {{\big( Y_- \, X \, 
{Y_-}^{\hskip-5pt -1} \big)}_+}^{\hskip-5pt +1} {\Big)}_-  =  \; 
\pm \pmatrix  \zeta^{-1} \Theta^{-1}  \quad  &  0  \\ \eta \, 
z^{+2} \Theta^{-1}  &  \zeta^{+1} \Theta^{+1}  \endpmatrix \cr }   
\eqno (4.7)  $$   
which gives  
  $$  \displaylines{ 
   J^{-1} \Big( {{\big( Y_- \, X \, {Y_-}^{\hskip-5pt -1}
\big)}_+}^{\hskip-5pt -1} \cdot Y \cdot {{\big( Y_- \, X \, 
{Y_-}^{\hskip-5pt -1} \big)}_+}^{\hskip-5pt +1} \Big) =   \hfill  
\cr 
   \hfill   = \pm \hskip-3pt \left( \hskip-3pt  \pmatrix  \zeta^{+1}
\Theta^{+1}  \;  &  \chi \, z^{-2} \Theta^{+1} + x \, z^{-1} 
\zeta^{+1} \Theta^{+1} - x \, z^{-1} \zeta^{-3} \Theta^{-1}  \\ 0  
&  \zeta^{-1} \Theta^{-1}  \endpmatrix  \hskip-2pt , \hskip-2pt 
\pmatrix  \zeta^{-1} \Theta^{-1}  \quad  &  0  \\  \eta \, z^{+2} 
\Theta^{-1}  &  \zeta^{+1} \Theta^{+1}  \endpmatrix  \hskip-3pt 
\right)  \cr }  $$   
as possible preimages of  $ \, {{\big( Y_- \, X \,
{Y_-}^{\hskip-5pt -1} \big)}_+}^{\hskip-5pt -1} \cdot Y \cdot 
{{\big( Y_- \, X \, {Y_-}^{\hskip-5pt -1} \big)}_+}^{\hskip-5pt 
+1} \, $  in  $ \, G^* \hskip-1,5pt \times G^* \, $.  This takes 
care of the second entry in the right-hand-side of (4.4).  
Finally, imposing the condition  $ \; \R_{{}_{W\!X}} 
\big(e_{G^*},e_{G^*}\big) = \big(e_{G^*},e_{G^*}\big) \; $  we 
must always take the "plus" signs, here and in (4.6).   
                                              \par   
   Using notation  
  $ \; x_1 := x \otimes 1 \, $,  $ \, {z_1}^{\pm 1} := z^{\pm 1}
\otimes 1 \, $,  $ \, y_1 := y \otimes 1 \, $,  $ \, x_2 := 1 
\otimes x \, $,  $ \, {z_2}^{\pm 1} := 1 \otimes z^{\pm 1} \; $ 
and  $ \; y_2 := 1 \otimes y \; $  we see that these last 
formul\ae{}  together with (4.6) give   
% 
% 
% 
%   $$  x_1 := x \otimes 1 \, ,  \quad  {z_1}^{\pm 1} := z^{\pm 1}
% \otimes 1 \, ,  \quad  y_1 := y \otimes 1 \, ,  \quad  x_2 := 1 
% \otimes x \, ,  \quad  {z_2}^{\pm 1} := 1 \otimes z^{\pm 1} \, , 
% \quad  y_2 := 1 \otimes y \, ,  $$   
% we see that these last formul\ae{}  together with (4.6) give   
% 
% 
% 
  $$  \hbox{ $ \eqalign{ 
   \R_{{}_{W\!X}} (x_1) = 
x_1 \cdot {z_2}^{\! -2} \cdot \Theta^{-1} \; ,  \qquad \qquad 
\R_{{}_{W\!X}} \big( {z_1}^{\!\pm 1} \big) =  {z_1}^{\!\pm 1} 
\cdot \Theta^{\mp 1}  \qquad  \cr  
   \R_{{}_{W\!X}} (y_1) = 
y_1 \cdot {z_2}^{\!+2} \cdot \Theta^{+1}  + y_2 \cdot {z_2}^{\!+1} 
{z_1}^{\!-1} \cdot \Theta^{+1} - y_2 \cdot {z_2}^{\!+1} 
{z_1}^{\!+3} \cdot \Theta^{-1}  \cr   
   \R_{{}_{W\!X}} (x_2) = 
x_2 \cdot {z_1}^{\!-2} \cdot \Theta^{+1}  + x_1 \cdot {z_1}^{\!-1} 
{z_2}^{\!+1} \cdot \Theta^{+1} - x_1 \cdot {z_1}^{\!-1} 
{z_2}^{\!-3} \cdot \Theta^{-1}  \cr   
   \R_{{}_{W\!X}} \big( {z_2}^{\!\pm 1} \big) = 
{z_2}^{\!\pm 1} \cdot \Theta^{\pm 1} \; ,  \qquad \qquad 
\R_{{}_{W\!X}} (y_2) =  y_2 \cdot {z_1}^{\! +2} \cdot \Theta^{-1} 
\qquad  \cr } $ }   \eqno (4.8)  $$   for the braiding  $ \, 
\R_{{}_{W\!X}} \, $.  To summarize, our discussion lead to the 
following result (which somewhat improves the analysis of the like 
problem performed in [WX], \S 9.7):   

\vskip9pt   

\proclaim{Theorem 4.4} \, Let  $ \, {\big( {}_s{{SL}_2}^{\!*} 
\hskip-1,5pt \times {}_s{{SL}_2}^{\!*} \big)}_2 \, $  be a
twofold covering of  $ \, {}_s{{SL}_2}^{\!*} \hskip-1,5pt \times 
{}_s{{SL}_2}^{\!*} \, $  and let  $  \, {\big( {}_s{{SL}_2}^{\!*} 
\hskip-1,5pt \!\times\! {}_s{{SL}_2}^{\!*} \big)}_2^{(\Theta)} 
\!\! := \! {\big( {}_s{{SL}_2}^{\!*} \hskip-1,5pt \!\times\! 
{}_s{{SL}_2}^{\!*} \big)}_2 \setminus \big\{ \Theta \not= 0 \big\} 
\, $  (a Zarisky open subset of  $ \, {}_s{{SL}_2}^{\!*} 
\hskip-1,5pt \!\times\! 
                   {}_s{{SL}_2}^{\!*} $).\hfill\break    
   \indent   Then the map  $ \Cal{R}_{{}_{W\!X}} $  is a Poisson
diffeomorphism from  $ \, {\big( {}_s{{SL}_2}^{\!*} \hskip-1,5pt 
\times {}_s{{SL}_2}^{\!*} \big)}_2^{(\Theta)} \, $  to itself.  
                                         \hfill\break    
% 
% Ou bien, \`a vrai
%                       dire, il suffirait de prendre la variet\'e 
%                       affine dont l'alg\`ebre de fonctions
%                       reguli\`eres soit  $ \; F \big[ {\big( 
%                       {}_s{{SL}_2}^{\!*} \hskip-1,5pt \times
%                       {}_s{{SL}_2}^{\!*}\big)}_2^{(\Theta)} \big] =
%                       {F \big[ {\big( {}_s{{SL}_2}^{\!*} \hskip-1,5pt 
%                       \times {}_s{{SL}_2}^{\!*} \big)}_2^{(\Theta)}
%                       \big]}_{\Theta} \; $  (localisation \`a 
%                       $ \Theta $)
% 
% 
% 
   \indent   In addition,  $ \Cal{R}_{{}_{W\!X}} $  is well
defined also on a distinguished variety  $ \; {\big(
{}_a{{SL}_2}^{\!*} \hskip-1,5pt \times {}_a{{SL}_2}^{\!*}
\big)}_2^{(\Theta)} \, $  which is a twofold covering of 
$ \, {}_a{{SL}_2}^{\!*} \hskip-1,5pt \times {}_a{{SL}_2}^{\!*}
\, $  minus one distinguished divisor.  In terms of function
algebras, these diffeomorphisms are uniquely determined by
formul\ae{} (4.8), which also define the braiding 
$ \; \R_{{}_{W\!X}} \, \colon \, F\big[\big[ {\gersl_2}^{\!*} 
\oplus {\gersl_2}^{\!*} \big]\big] \;{\buildrel \cong \over 
{\relbar\joinrel\relbar\joinrel\llongrightarrow}}\; F \big[\big[ 
{\gersl_2}^{\!*} \oplus {\gersl_2}^{\!*} \big]\big] \; $.   \qed   
\endproclaim   

\vskip9pt 

  {\bf 4.5 The quantisation deformation construction.} {\it
(Warning: in the present section we follow the lines of [Ga1], but 
we adopt  {\sl different normalisations}  in the definition of 
quantum groups and their  $ R $--matrices)}  \, Let  $ \, 
U_\hbar(\gerg) = U_\hbar(\gersl_2) \, $  be the unital associative 
topological  $ \kh $--algebra  with (topological) generators $ \, 
X $,  $ H $,  $ Y $,  and relations    
  $$  H X - X H = + 2 \, X \, ,  \quad  H Y - Y H = - 2 \, Y \, , 
\quad  X Y - Y X = {\, e^{+\hbar H / 2} - e^{-\hbar H / 2} \, 
\over \, e^{+\hbar / 2} - e^{-\hbar / 2} \,} \; .   \eqno (4.9)  
$$   For later use we set also  $ \; L^{\pm 1} := e^{\pm \hbar H / 
4} \; $  and  $ \; q^{\pm 1} := e^{\pm \hbar / 2} \, $;  \; 
therefore  
  $$  L^{\pm 1} X = q^{\pm 1} X L^{\pm 1} \, ,  \quad  L^{\pm 1} Y
= q^{\mp 1} Y L^{\pm 1} \, ,  \quad  X Y - Y X = {\, L^{+2} - 
L^{-2} \, \over \, q^{+1} - q^{-1} \,} \; .   \eqno (4.10)  $$    
   \indent   There is a Hopf algebra structure on 
$ U_\hbar(\gersl_2) $,  given on generators by   
  $$  \displaylines{ 
   \Delta(X) = X \otimes e^{+\hbar H / 4} + e^{-\hbar H / 4}
\otimes X = X \otimes L^{+1} + L^{-1} \otimes X  \cr  
   \Delta(H) = H \otimes 1 + 1 \otimes H  \cr  
   \Delta(Y) = Y \otimes e^{+\hbar H / 4} + e^{-\hbar H / 4}
\otimes Y = Y \otimes L^{+1} + L^{-1} \otimes Y  \cr  
   \epsilon(X) = \epsilon(H) = \epsilon(Y) = 0 \, ,  \quad 
\epsilon\big(L^{\pm 1}\big) = 1  \cr   
   S(X) = - e^{-\hbar/2} X = - q^{-1} X \, ,  \; S(H) = -H \, ,  \;
S(Y) = - e^{-\hbar/2} Y = - q^{-1} Y \, ,  \; S\big(L^{\pm 1}\big) 
= L^{\mp 1} .  \cr }    $$  Then  $ U_\hbar(\gersl_2) $  is a 
QUEA, whose semiclassical limit is $ U(\gersl_2) $  (w.r.t.~the 
co-Poisson structure considered in \S 4.2). For later use we 
record that  
  $$  \Big\{\, X^a H^b Y^c \,\Big\vert\, a, b, c \in \N \,\Big\}
\hskip6pt  \hbox{\sl is a topological  $ \kh $--basis  of} 
\hskip5pt  U_\hbar(\gersl_2) \, .   \eqno (4.11)  $$   
   \indent   The very definitions also show that the unital subalgebra
of  $ U_\hbar(\gersl_2) $  generated over the Laurent polynomial 
ring $ \kqqm $  by  $ \; X \, $,  $ \, L^{\pm 1} \, $,  $ \, D := 
(L-1) \Big/ (q-1) \, $,  $ \; \varGamma := \big( L^{+2} - L^{-2} 
\big) \Big/ \big( q^{+1} - q^{-1} \big) \; $  and  $ \, Y \; $ 
is a Hopf algebra (over $ \kqqm $)  as well, which we denote by 
$ \, U_q^s(\gersl_2) \, $.  Similarly, the unital subalgebra of 
$ U_\hbar(\gersl_2) $ generated over the Laurent polynomial ring 
$ \kqqm $  by  $ \; X L^{-1} \, $, $ \, K^{\pm 1} := L^{\pm 2}
\, $,  $ \; T := (K-1) \Big/ (q-1) \, $,  $ \; \varGamma := \big(
K^{+1} - K^{-1} \big) \Big/ \big( q^{+1} - q^{-1} \big) \; $ 
and  $ \, L^{+1} Y \; $  is a Hopf algebra as well (a Hopf 
$ \kqqm $--subalgebra  of  $ U_q^s(\gersl_2) \, $),  which
we denote by $ \, U_q^a(\gersl_2) \, $.   
                                            \par   
   Now we go and compute  $ \, {U_\hbar(\gersl_2)}' $.  From
definitions we get, for any  $ \, n \in \N $,  
  $$  \displaylines{ 
   {} \quad  \delta_n(X) = {(\hbox{\rm id} - \epsilon)}^{\otimes n}
\big( \Delta^n(X) \big) = {(\hbox{\rm id} - \epsilon)}^{\otimes n}
\left( \sum_{s=1}^n {\big( L^{-1} \big)}^{\otimes (s-1)} \otimes
X \otimes {\big( L^{+1} \big)}^{\otimes (n-s)} \right) =  
\hfill  \cr  
   {} \hfill   = \sum_{s=1}^n {\left( \sum_{t>0} {(-\hbar)}^t H^t
\Big/ t! \right)}^{\hskip-3pt \otimes (s-1)} \hskip-5pt \otimes
X \otimes {\left( \sum_{r>0} {(+\hbar)}^r H^r \Big/ r! 
\right)}^{\hskip-3pt \otimes (n-s)}  \hskip-9pt  \in \hskip5pt 
\hbar^{n-1} U_\hbar(\gersl_2) \setminus \hbar^n U_\hbar(\gersl_2)  
\cr }  $$   from which we get  $ \, \dot{X} := \hbar \, X \in 
{U_\hbar(\gerg)}' \setminus \hbar \, {U_\hbar(\gerg)}' \, $.  
Similarly  $ \, \dot{Y} := \hbar \, Y \in {U_\hbar(\gerg)}' 
\setminus \hbar \, {U_\hbar(\gerg)}' \, $.  As for the generator  
$ H $,  we have $ \; \Delta^n(H) = \sum_{s=1}^n 1^{\otimes (s-1)} 
\otimes H \otimes 1^{\otimes (n-s)} \; $  for all  $ \, n \in
\N \, $,  \; whence for $ \, \delta_n = {(\hbox{\rm id} -
\epsilon)}^{\otimes n} \circ \Delta^n \, $  we have   
  $$  \delta_0(H) = 0 \, ,  \qquad  \delta_1(H) = H \in
U_\hbar(\gerg) \setminus \hbar \, U_\hbar(\gerg) \, ,  \qquad 
\delta_n(H) = 0 \in \hbar^n U_\hbar(\gersl_2)  \quad\;  \forall
\;\, n > 1 \, ,  $$   
so that  $ \, \dot{H} := \hbar \, H \in {U_\hbar(\gersl_2)}' 
\setminus \hbar \, {U_\hbar(\gersl_2)}' \, $.  Therefore 
$ {U_\hbar(\gersl_2)}' $  contains the subalgebra  $ U' $ 
topologically generated by  $ \, \dot{X} $,  $ \dot{H} $, 
$ \dot{Y} \, $.  On the other hand, using (4.11) a thorough  
--- but straightforward ---   computation shows that any element
in  $ {U_\hbar(\gersl_2)}' $  does necessarily lie in  $ U' $ 
(details are left to the reader: everything follows from 
definitions and the formulas for  $ \, \Delta^n \, $).  Thus $ 
{U_\hbar(\gersl_2)}' $  is nothing but the unital subalgebra of  $ 
U_\hbar(\gersl_2) $  topologically generated by  $ \, \dot{X} $,  
$ \dot{H} $,  $ \dot{Y} \, $.  As a consequence, $ 
{U_\hbar(\gersl_2)}' $  can be presented as the unital associative 
topological  $ \kh $--algebra  with (topological) generators  $ \, 
\dot{X} $,  $ \dot{H} $,  $ \dot{Y} \, $ and relations  
  $$  \hbox{ $ \eqalign{ 
   \dot{H} \dot{X} - \dot{X} \dot{H} = + 2 \hbar \, \dot{X}
\, ,  \qquad  \dot{H} \dot{Y}  &  - \dot{Y} \dot{H} = - 2 \hbar \, 
\dot{Y}  \cr   
   \dot{X} \dot{Y} - \dot{Y} \dot{X} \, = \, \hbar \, A \cdot
\Big( e^{+\dot{H}/2} - e^{-\dot{H}/2} \Big)  &  = \, \hbar \, A 
\cdot \big( L^{+2} - L^{-2} \big)  \cr } $ }   \eqno (4.12)  $$   
where  $ \; A := \hbar^2 \Big/ \big( e^{+\hbar/2} - e^{-\hbar/2} 
\big) = \hbar \cdot {\Big( \sum_{s > 0} {(+\hbar/2)}^{2s} \Big/ 
(2s-1)! \Big)}^{\!-1} \; \big( \in \kh \, \big) \, $,  \; with 
Hopf algebra structure given by   
  $$  \displaylines{ 
   \Delta\big(\dot{X}\big) = \dot{X} \otimes e^{+\dot{H}/4}
+ e^{-\dot{H}/4} \otimes \dot{X} = \dot{X} \otimes L^{+1} + L^{-1} 
\otimes \dot{X}  \cr   
   \Delta\big(\dot{H}\big) =
\dot{H} \otimes 1 + 1 \otimes \dot{H}  \cr  
   \Delta\big(\dot{Y}\big) = \dot{Y} \otimes e^{+\dot{H}/4}
+ e^{-\dot{H}/4} \otimes \dot{Y} = \dot{Y} \otimes L^{+1} + L^{-1} 
\otimes \dot{Y}  \cr   
   \epsilon\big(\dot{X}\big) = \epsilon\big(\dot{H}\big)
= \epsilon\big(\dot{Y}\big) = 0 \, ,  \quad \epsilon\big(L^{\pm 
1}\big) = 1  \cr   
   S\big(\dot{X}\big) = - e^{-\hbar/2} \dot{X} = - q^{-1} \dot{X} , 
\; S\big(\dot{H}\big) = -\dot{H} ,  \; S\big(\dot{Y}\big) = - 
e^{-\hbar/2} \dot{Y} = - q^{-1} \dot{Y} ,  \; S\big(L^{\pm 1}\big) 
= L^{\mp 1} .  \cr }    $$  
   \indent   As an immediate consequence, this description yields
also a similar presentation of  $ \; {U_\hbar(\gersl_2)}' \Big/ 
\hbar \, {U_\hbar(\gersl_2)}' \, $:  \; then comparing the latter 
with the presentation of  $ F \big[ \big[{\gersl_2}^{\!*}\big]
\big] $  that one argues from \S 4.2 we find that, as predicted
by the quantum duality principle (cf.~Theorem 1.6)  {\sl there
is an isomorphism of (topological) Poisson Hopf algebras}   
  $$  \varPhi_h \, \colon \, {U_\hbar(\gersl_2)}' \Big/ \hbar \,
{U_\hbar(\gersl_2)}' = \Bbbk \Big[ \Big[ \dot{X}{\big|}_{\hbar=0}
\, , \dot{H}{\big|}_{\hbar=0} \, , \dot{Y}{\big|}_{\hbar=0} \Big]
\Big] \; {\buildrel \cong \over {\llongrightarrow}} \; F\big[\big[
{\gersl_2}^{\!*} \big]\big] = k\big[\big[x,(z-1),y\big]\big]  $$  
where we set  $ \; S{\big|}_{\hbar=0} := S \mod \hbar \, 
{U_\hbar(\gersl_2)}' \, $  for all  $ \, S \in 
{U_\hbar(\gersl_2)}' \; $  and the Poisson structure considered on  
$ \, {U_\hbar(\gersl_2)}' \Big/ \hbar \, {U_q^s(\gersl_2)}' \, $  
is the one given by the standard recipe (see  \S 1.3{\it (b)})   
  $$  \big\{ a{\big|}_{\hbar=0} \, , b{\big|}_{\hbar=0} \big\}
:= {\left( {\, a \, b - b \, a \, \over \, \hbar \,} \right)} 
{\bigg|}_{\hbar=0}  \qquad  \forall \quad a, b \in 
{U_\hbar(\gersl_2)}' \; ;  $$   
explicitly,  $ \varPhi_\hbar $  is given by   
  $$  \dot{X}{\big|}_{\hbar=0} \mapsto x \; ,  \qquad 
\dot{H}{\big|}_{\hbar=0} \mapsto -4 \hskip0,1pt \log(z) \; , 
\qquad  L^{\pm 1}{\big|}_{\hbar=0} \mapsto z^{\mp 1} \; , \qquad  
\dot{Y}{\big|}_{\hbar=0} \mapsto y \; .   \eqno (4.13)  $$   
   \indent   Note also that the unital  $ \kqqm $--subalgebra 
of  $ U_\hbar(\gersl_2) $   --- and of  $ U_q^s(\gersl_2) $  ---  
generated by  $ \; \check{X} := (q-1) \, X \, $,  $ \, L^{\pm 1} 
\, $,  $ \; \check{\varGamma} := (q-1) \, \varGamma \; $  and $ \, 
\check{Y} := (q-1) \, Y \; $  is in fact a  {\sl Hopf\/} 
subalgebra, which we denote by  $ \, {U_q^s(\gersl_2)}' \, $ (note 
also that  $ \; \check{D} := (q-1) \, D = L - 1 \in 
{U_q^s(\gersl_2)}' \, $  too).  Indeed,  $ \, {U_q^s(\gersl_2)}' 
\, $  admits the presentation by the above generators and 
relations  
  $$  \displaylines{
   L^{\pm 1} L^{\mp 1} = 1 \, ,  \,\;  L^{\pm 1} \check{\varGamma}
= \check{\varGamma} L^{\pm 1} ,  \,\;  \big( 1 + q^{-1} \big) 
\check{\varGamma} = L^{+2} - L^{-2} ,  \,\;  \check{X} \check{Y} - 
\check{Y} \check{X} = (q-1) \check{\varGamma}  \cr   
   L^{+2} - L^{-2} = \big( 1 + q^{-1}\big) \check{\varGamma} \, , 
\quad  L^{\pm 1} \check{Y} = q^{\mp 1} \check{Y} L^{\pm 1} \, , 
\quad  L^{\pm 1} \check{X} = q^{\pm 1} \check{X} L^{\pm 1}  \cr  
   \check{\varGamma} \check{Y} = q^{-2} \check{Y} \check{\varGamma}
- (q-1) \big( q + q^{-1} \big) \check{F} \, ,  \qquad 
\check{\varGamma} \check{X} = q^{+2} \check{X} \check{\varGamma} + 
(q-1) \big( q + q^{-1} \big) \check{X}  \cr }  $$   with Hopf 
structure given by   
  $$  \matrix
      \Delta\big(\check{X}\big) = \check{X} \otimes L^{+1} + L^{-1}
\otimes \check{X} \, ,  &  \hskip25pt  \epsilon\big(\check{X}\big) 
= 0 \, ,  &  \hskip25pt  S\big(\check{X}\big) = - q^{-1} \check{X} 
\, \phantom{.}  \\   
      \Delta\big(\check{\varGamma}\big) = \check{\varGamma} \otimes
L^{+2} + L^{-2} \otimes \check{\varGamma} \, ,  &  \hskip25pt 
\epsilon\big(\check{\varGamma}\big) = 0 \, ,  &  \hskip25pt 
S\big(\check{\varGamma}\big) = - \check{\varGamma}  \\   
      \Delta \big( L^{\pm 1} \big) = L^{\pm 1} \otimes L^{\pm 1} \, , 
&  \hskip25pt  \epsilon \big( L^{\pm 1} \big) = 1 \, ,  &  
\hskip25pt S \big( L^{\pm 1} \big) = L^{\mp 1} \, \phantom{.}  \\   
      \Delta\big(\check{Y}\big) = \check{Y} \otimes L^{+1} + L^{-1}
\otimes \check{Y} \, ,  &  \hskip25pt  \epsilon\big(\check{Y}\big) 
= 0 \, ,  &  \hskip25pt  S\big(\check{Y}\big) = - q^{-1} \check{Y} 
\, .  \\   
     \endmatrix  $$
   \indent   Similarly, the unital  $ \kqqm $--subalgebra  of 
$ U_\hbar(\gersl_2) $   --- and of  $ U_q^a(\gersl_2) $  and $ 
{U_q^s(\gersl_2)}' $  ---   generated by  $ \; \check{X} L^{-1} \, 
$,  $ \, K^{\pm 1} \, $,  $ \; \check{\varGamma} \; $  and $ \, 
L^{+1} \check{Y} \; $  is in fact a  {\sl Hopf\/}  subalgebra too, 
which we denote by  $ \, {U_q^a(\gersl_2)}' \, $,  \, and is of 
course a Hopf subalgebra of  $ \, {U_q^s(\gersl_2)}' \, $ as well 
(with  $ \; \check{T} := (q-1) \, T = K - 1 \in {U_q^a(\gersl_2)}' 
\; $  too).   
                                            \par  
   Now this description yields also a similar presentation of 
$ \; {U_q^s(\gersl_2)}' \Big/ (q-1) \, {U_q^s(\gersl_2)}' \, $: \; 
then comparing the latter with the presentation of  $ F \big[ 
{}_s{{SL}_2}^{\!*}\big] $  in \S 4.2 we find that  {\sl there is a 
Poisson Hopf algebra isomorphism}    
  $$  \varPhi_q^s \, \colon \, {U_q^s(\gersl_2)}' \Big/ (q-1) \,
{U_q^s(\gersl_2)}' \; {\buildrel \cong \over {\relbar\joinrel 
\relbar\joinrel\relbar\joinrel\relbar\joinrel\relbar 
\joinrel\relbar\joinrel\llongrightarrow}} \; F\big[ 
{}_s{{SL}_2}^{\!*}\big]  $$   where we set  $ \; S{\big|}_{q=1} := 
S \mod (q-1) \, {U_q^s(\gersl_2)}' \, $  for all  $ \, S \in 
{U_q^s(\gersl_2)}' \; $ and the Poisson structure considered on  $ 
\, {U_q^s(\gersl_2)}' \Big/ (q-1) \, {U_q^s(\gersl_2)}' \, $  is 
the one given by the standard recipe   
  $$  \big\{ a{\big|}_{q=1} \, , b{\big|}_{q=1} \big\} \; := \;
{\left( {\, a \, b - b \, a \, \over \, (q-1) \,} \right)}
{\bigg|}_{q=1} \qquad  \forall \quad a, b \in {U_q^s(\gersl_2)}'
\; ;  $$   
explicitly,  $ \varPhi_q^s $  is given by   
  $$  \check{X}{\big|}_{q=1} \mapsto x \big/ 2 \; ,  \qquad 
\check{\varGamma}{\big|}_{q=1} \mapsto \big( z^{-2} - z^{+2} \big) 
\Big/ 2 \; ,  \qquad  L^{\pm 1}{\big|}_{q=1} \mapsto z^{\mp 1} \; 
, \qquad  \check{Y}{\big|}_{q=1} \mapsto y \big/ 2 \; .  $$   
  \indent   In addition,  $ \varPhi_q^s $  gives by restriction
a similar Poisson Hopf algebra isomorphism   
  $$  \displaylines{ 
   \varPhi_q^a \, \colon \, {U_q^a(\gersl_2)}' \Big/ (q-1) \,
{U_q^a(\gersl_2)}' \; {\buildrel \cong \over {\relbar\joinrel 
\relbar\joinrel\relbar\joinrel\relbar\joinrel\relbar\joinrel\relbar 
\joinrel\llongrightarrow}} \; F\big[{}_a{{SL}_2}^{\!*}\big]  \cr   
   {\big( \check{X} L^{-1} \big)}{\big|}_{q=1} \! \mapsto x z \big/ 2
\, ,  \; \check{\varGamma}{\big|}_{q=1} \! \mapsto \big( z^{-2} \! 
- \! z^{+2} \big) \Big/ 2 \, ,  \; K^{\pm 1}{\big|}_{q=1} \! 
\mapsto z^{\mp 2} ,  \; {\big( L^{+1} \check{Y} 
\big)}{\big|}_{q=1} \! \mapsto z^{-1} y \big/ 2 \, .  \cr }  $$   
   \indent   The reason for considering  $ \, U_q^c(\gersl_2) \, $ 
and  $ \, {U_q^c(\gersl_2)}' $  (for  $ \, c = a, s \, $)  is that 
we can compute the braiding  $ \, \R_{{}_{G\!H}} \, $  through 
them, as we shall see in the sequel.   
                                            \par  
   First,  $ U_\hbar(\gersl_2) $  is indeed a \QTQUEA, whose 
$ R $--matrix  is  $ \; R_\hbar = R_0 \cdot R_1 \; $  with  
  $$  R_0 = \exp \Big( \hbar \cdot H \otimes H \big/ 4 \Big) \, , 
\hskip13pt  R_1 = \sum_{n \in \N} {\; {\big( e^\hbar \big)}^{{n+1} 
\choose 2} \, \over \; {(n)}_{e^\hbar}! \,} \, {\big( e^\hbar - 1 
\big)}^n \cdot {\big( e^{+\hbar H / 4} X \big)}^{\!n} \! \otimes 
{\big( e^{-\hbar H / 4} Y \big)}^{\!n}  $$   where  $ \, {(n)}_a! 
:= \prod_{r=1}^n {\, a^n - 1 \, \over \, a - 1 \,} \, $  (in this 
case  $ \, a = {e^\hbar} \, $).  This $ R $-matrix  is a 
quantisation of the classical  $ \r $--matrix of  $ \gersl_2 $,  
in the sense that  $ \, R_\hbar = 1 + \r \, \hbar + \O \left( 
\hbar^2 \right) \, $,  \, where  $ \, \O \left( \hbar^2 \right)
\, $  is some element of  $ \, \hbar^2 \cdot \uhg \otimes \uhg
\, $  (like in Remark 1.9{\it (b)\/});  thus the \QTQUEA{}  $ \,
\big( U_\hbar(\gersl_2), R_\hbar \big) \, $  is a quantisation of
the quasitriangular Lie bialgebra  $ \, (\gersl_2, \r \big) \, $, 
\, as required to ignite the quantisation deformation procedure
to construct a braiding on  $ F[[\gerg^* \oplus \gerg^*]] $ 
for  $ \, \gerg = \gersl_2 \, $.   
                                            \par  
   Now, we are interested in the braiding operator  $ \,
\R_{{}_{G\!H}} \, $  induced at  $ \, \hbar = 0 \, $  by the 
operator  $ \, \R_\hbar := {\hbox{\rm Ad}}(R_\hbar) \, $  acting 
on the algebra  $ {\big( U_\hbar(\gersl_2) \otimeshat U_\hbar
(\gersl_2) \big)}' \! = {U_\hbar(\gersl_2)}' \otimestilde
{U_\hbar(\gersl_2)}' \, $.     
                                            \par  
   We perform the calculation along the following lines.  As the 
$ R $--matrix  factors into  $ \, R_\hbar = R_0 \cdot R_1 \, $, 
\, we compute separately the adjoint action of the two
factors onto  $ \, {U_\hbar(\gersl_2)}' \otimestilde
{U_\hbar(\gersl_2)}' $  modulo  $ \hbar $.  A first analysis
shows that both actions are given by exponentials of Hamiltonian
vector fields on the formal Poisson group  $ {\gersl_2}^{\!*}
\times {\gersl_2}^{\!*} $.  The first action   --- namely,
that arising from  $ R_0 $  ---   is computed via
straightforward calculation.  As for the second action  
--- the one of  $ R_1 $  ---   one in fact has to compute
the action of a Hamiltonian vector field on  $ \,
{}_s{{SL}_2}^{\!*} \times {}_s{{SL}_2}^{\!*} \, $ 
(minus a divisor): using Leibniz' rule, one reduces
to compute the action of some Hamiltonian vector
fields on  $ \, {}_s{{SL}_2}^{\!*} \, $  alone.   
                                           \par   
   To begin with, write  $ \, R_\hbar = R_0 \cdot R_1 \, $  in
terms of  $ {U_\hbar(\gersl_2)}' $;  like in [Ga1], \S 3, we find  
  $$  \displaylines{ 
   R_0 = \exp \Big( \hbar \cdot H \otimes H \big/ 4 \Big) = \exp
\Big( \hbar^{-1} \cdot \dot{H} \otimes \dot{H} \big/ 4 \Big)  \cr   
   R_1 = \sum_{n \in \N} {\; {\big( e^\hbar \big)}^{{n+1} \choose 2}
\over \; {(n)}_{e^\hbar}!} \, {\big( e^\hbar \hskip-1pt - 
\hskip-1pt 1 \big)}^n \cdot {\big( e^{+\hbar H / 4} X \big)}^{\!n} 
\! \otimes {\big( e^{-\hbar H / 4} Y \big)}^{\!n} \hskip-2pt =  
{\left( \hskip-2pt {\big( e^\hbar \hskip-1pt - \hskip-1pt 1 
\big)}^2 \hskip-3pt \cdot \hskip-2pt L^{+1} X \hskip-2pt \otimes 
\hskip-2pt L^{-1} Y \, ; \, e^\hbar \right)}_\infty  \cr }  $$   
where  $ \; {(z;q)}_\infty := \prod_{n \in \N} (1 - z \, q^n) \, 
$. Now, the behaviour of  $ R_1 $  when  $ \, \hbar \rightarrow 0 
\, $ is ruled by [Re], Lemma 3.4.1 (see also [Ga1], Lemma 2.2): 
namely (proceeding as in [Ga1], \S 3), we have   
  $$  \eqalign{
   R_1  &  \; = \; \exp \left({\, -1 \, \over \, \hbar \,} \cdot
{\int_0}^{{( e^\hbar - 1 )}^2 \cdot L^{+1} X \otimes L^{-1} Y} {\, 
\log(1-\tau) \, \over \, \tau \,} \; d\tau \; \cdot \big( 1 + 
\hbar \, C \big) \right) =  \cr   
   &  \; = \; \exp \left({\, -1 \, \over \, \hbar \,} \cdot
{\int_0}^{L^{+1} \dot{X} \otimes L^{-1} \dot{Y}} {\, \log(1+t) \, 
\over \, t \,} \; dt \; \cdot \big( 1 + \hbar \, K \big) \right) 
\cr }  $$   where  $ \; {\int_0}^{L^{+1} \dot{X} \otimes L^{-1} 
\dot{Y}} {\, \log(1+t) \, \over \, t \,} \; dt := \sum_{n>0} \big( 
L^{+1} \dot{X} \otimes L^{-1} \dot{Y} \big) \Big/ n^2 \; $ (use 
Mac Laurin expansion of  $ \, \log(1+x) \, $)  and  $ C $ and  $ K 
$  denote a suitable elements of  $ \, {U_\hbar(\gersl_2)}' 
\otimestilde {U_\hbar(\gersl_2)}' \, $,  \, namely again power 
series in  $ \, L^{+1} \dot{X} \otimes L^{-1} \dot{Y} \, $,  \, 
hence they commute with  $ \, {\int_0}^{L^{+1} \dot{X} \otimes 
L^{-1} \dot{Y}} {\, \log(1-t) \, \over \, t \,} \; dt \; $;  \, so   
  $$  \eqalign{ 
   R_1  &  \; = \; \exp \left({\, -1 \, \over \, \hbar \,} \cdot
{\int_0}^{L^{+1} \dot{X} \otimes L^{-1} \dot{Y}} {\, \log(1+t) \, 
\over \, t \,} \; dt \; \cdot \big( 1 + \hbar \, K \big) \right) =  
\cr   
   &  \; = \; \exp \left({\, -1 \, \over \, \hbar \,} \cdot
{\int_0}^{L^{+1} \dot{X} \otimes L^{-1} \dot{Y}} {\, \log(1+t) \, 
\over \, t \,} \; dt \right) \cdot Z  \cr }  $$   for some  $ \; Z 
\in {U_\hbar(\gersl_2)}' \otimestilde {U_\hbar(\gersl_2)}' \, $.  
Of course we have   
  $$  \R_\hbar := {\hbox{\rm Ad}}(R_\hbar) = {\hbox{\rm Ad}}
(R_0 \cdot R_1) = {\hbox{\rm Ad}}(R_0) \circ {\hbox{\rm Ad}}(R_1) 
= \R_\hbar^{(0)} \circ \R_\hbar^{(1)}  $$   with  $ \; 
\R_\hbar^{(0)} := {\hbox{\rm Ad}}(R_0) \, $, $ \; \R_\hbar^{(1)} 
:= {\hbox{\rm Ad}}(R_1) \, $.  Thus also   
  $$  \R_{{}_{G\!H}} := \R_\hbar{\Big|}_{\hbar=0} =
\R_{{}_{G\!H}}^{(0)} \circ \R_{{}_{G\!H}}^{(1)}  \qquad  
\hbox{with} \quad  \R_{{}_{G\!H}}^{(0)} := 
\R_\hbar^{(0)}{\Big|}_{\hbar=0} \, , \quad  \R_{{}_{G\!H}}^{(1)} 
:= \R_\hbar^{(1)}{\Big|}_{\hbar=0} \; .  \eqno (4.14)  $$    
Finally we have   
  $$  \eqalign{ 
   \R_\hbar^{(1)}  &  := {\hbox{\rm Ad}}(R_1) = {\hbox{\rm Ad}}
\Bigg( \exp \left({\, -1 \, \over \, \hbar \,} \cdot 
{\int_0}^{L^{+1} \dot{X} \otimes L^{-1} \dot{Y}} {\, \log(1+t) \, 
\over \, t \,} \; dt \right) \cdot Z \Bigg) =   \hfill  \cr   
   &  = {\hbox{\rm Ad}} \Bigg( \exp \left({\, -1 \, \over \, \hbar \,}
\cdot {\int_0}^{L^{+1} \dot{X} \otimes L^{-1} \dot{Y}} {\, 
\log(1+t) \, \over \, t \,} \; dt \right) \Bigg) \circ {\hbox{\rm 
Ad}}(Z) = \cr   
   &  = {\hbox{\rm Ad}} \Bigg( \exp \left({\, -1 \, \over \, \hbar
\,} \cdot {\int_0}^{L^{+1} \dot{X} \otimes L^{-1} \dot{Y}} {\, 
\log(1+t) \, \over \, t \,} \; dt \right) \Bigg)  \; \mod \; \hbar 
\cdot {U_\hbar(\gersl_2)}' \otimes {U_\hbar(\gersl_2)}'  \cr }  $$  
because  $ \; Z \in {U_\hbar(\gersl_2)}' \otimestilde 
{U_\hbar(\gersl_2)}' \, $  and  $ \, {\left( {U_\hbar(\gersl_2)}'
\otimestilde {U_\hbar(\gersl_2)}' \,\right)}{\!\Big|_{\hbar=0}} =
F\big[\big[ {\gersl_2}^{\!*} \times {\gersl_2}^{\!*} \big]\big] \, $ 
is commutative (hereafter, by  $ \, S{\big|}_{\hbar=0} \, $  we shall
denote the coset of  $ \, S \in {U_\hbar(\gersl_2)}' \otimestilde
{U_\hbar(\gersl_2)}' \, $  modulo  $ \, \hbar \cdot
{U_\hbar(\gersl_2)}' \otimestilde {U_\hbar(\gersl_2)}') \, $. 
Hence 
our analysis shows that  $ \, \R_\hbar^{(i)} = {\hbox{\rm Ad}} 
\Big( \exp \big( \hbar^{-1} \Lambda_i \big) \Big) \, $   with  $ 
\, \Lambda_i \in {U_\hbar(\gersl_2)}' \otimestilde 
{U_\hbar(\gersl_2)}' \, $  for $ \, i = 0, 1 \, $:  \, indeed, we 
found   
  $$  \Lambda_0 = \big( \dot{H} \otimes \dot{H} \big) \big/ 4 \; , 
\quad  \Lambda_1 = - {\int_0}^{L^{+1} \dot{X} \otimes L^{-1} 
\dot{Y}} {\, \log(1+t) \, \over \, t \,} \; dt \; = \, - 
\sum_{n>0} {\big( L^{+1} \dot{X} \otimes L^{-1} \dot{Y} \big)}^n 
\Big/ n^2 \; .  $$   
But then we have  $ \; \R_\hbar^{(i)} = {\hbox{\rm Ad}} \Big( \exp
\big( \hbar^{-1} \Lambda_i \big) \Big) = \exp \! \Big( {\hbox{\rm
ad}}_{[\ ,\ ]} \big( \hbar^{-1} \Lambda_i \big) \Big) = \exp \!
\Big( {\hbox{\rm ad}}_{{\, [\ ,\ ] \, \over \, \hbar \,} } \big( 
               \Lambda_i \big) \Big) \equiv $\break   
     $ \equiv \exp \! \Big( {\hbox{\rm ad}}_{\{\ ,\ \}} \big(
\Lambda_i{\big|}_{\hbar=0} \big) \Big) \, $,  \, that is  $ \; 
\R_{{}_{G\!H}}^{(i)} = \exp \! \Big( {\hbox{\rm ad}}_{\{\ ,\ \}} 
\big( \Lambda_i{\big|}_{\hbar=0} \big) \Big) \, $:  \, in geometric
terms, this means 
that  $ \R_{{}_{G\!H}}^{(i)} $  (hence also  $ \R_{{}_{G\!H}} $) 
is the integration of a Hamiltonian vector fields over the formal 
Poisson group  $ {\gersl_2}^{\!*} \times {\gersl_2}^{\!*} $.   
                                            \par  
   To describe  $ \R_{{}_{G\!H}}^{(0)} $  and 
$ \R_{{}_{G\!H}}^{(1)} $  we set  $ \, S_1 := S \otimes 1 \, $, $ 
\, S_2 := 1 \otimes S \, $  for any  $ \, S \in 
{U_\hbar(\gersl_2)} \, $  (note that  $ S_1 $  and  $ S_2 $  
commute with each other) and also  $ \overline{S} $  for any coset 
modulo  $ \hbar $.  
                                            \par  
   The case of  $ \R_{{}_{G\!H}}^{(0)} $  is trivial: direct
computation   --- using (4.12) ---  gives   
  $$  \eqalign{ 
   \R_\hbar^{(0)}\big(\dot{X}_1\big) = \dot{X}_1 \, L_2^{+2} \, , 
&  \hskip9pt  \R_\hbar^{(0)}\big(\dot{H}_1\big) = \dot{H}_1 \, , 
\hskip9pt  \R_\hbar^{(0)}\big(\dot{L}_1^{\pm 1}\big) = 
\dot{L}_1^{\pm 1} \, ,  \hskip9pt  
\R_\hbar^{(0)}\big(\dot{Y}_1\big) = \dot{Y}_1 \, L_2^{-2}  \cr   
   \R_\hbar^{(0)}\big(\dot{X}_2\big) = L_1^{+2} \, \dot{X}_2 \, , 
&  \hskip9pt  \R_\hbar^{(0)}\big(\dot{H}_2\big) = \dot{H}_2 \, , 
\hskip9pt  \R_\hbar^{(0)}\big(\dot{L}_2^{\pm 1}\big) = 
\dot{L}_2^{\pm 1} \, ,  \hskip9pt  
\R_\hbar^{(0)}\big(\dot{Y}_2\big) = L_1^{-2} \, \dot{Y}_2  \cr }  
$$    whence using (4.13) we argue at once for  $ \, 
\R_{{}_{G\!H}}^{(0)} \, \colon \, F \big[\big[ {\gersl_2}^{\!*} 
\oplus {\gersl_2}^{\!*} \big]\big] \;{\buildrel \cong \over 
{\llongrightarrow}}\; F \big[\big[ {\gersl_2}^{\!*} \oplus 
{\gersl_2}^{\!*} \big]\big] \, $  
  $$  \hbox{ $ \eqalign{ 
   \R_{{}_{G\!H}}^{(0)}(x_1) = x_1 \, z_2^{-2} \, ,  &  \hskip25pt 
\R_{{}_{G\!H}}^{(0)}\big(z_1^{\pm 1}\big) = z_1^{\pm 1} \, , 
\hskip25pt  \R_{{}_{G\!H}}^{(0)}(y_1) = y_1 \, z_2^{+2}  \cr   
   \R_{{}_{G\!H}}^{(0)}(x_2) = z_1^{-2} \, x_2 \, ,  &  \hskip25pt 
\R_{{}_{G\!H}}^{(0)}\big(z_2^{\pm 1}\big) = z_2^{\pm 1} \, , 
\hskip25pt  \R_{{}_{G\!H}}^{(0)}(y_2) = z_1^{+2} \, y_2  \cr } $ }  
\eqno (4.15)  $$    
(recall that  $ \; F \big[\big[ {\gersl_2}^{\!*} \oplus
{\gersl_2}^{\!*} \big]\big] = \Bbbk \big[\big[ x_1,
(z_1 - 1), y_1, x_2, (z_2 - 1), y_2 \big]\big] \; $  thus 
$ \R_{{}_{G\!H}}^{(0)} $  is uniquely determined by the
images of  $ x_i $,  $ z_i $,  $ y_i $  [$ i = 1, 2 $]).  
                                            \par  
   As for  $ \R_{{}_{G\!H}}^{(1)} $,  we proceed in steps.  First,
using the Jacobi identity for  $ \, \{\ ,\ \} \, $  we get   
  $$  \displaylines{ 
   \R_{{}_{G\!H}}^{(1)} = \exp \! \left( {\hbox{\rm ad}}_{\{\ ,\ \}}
\left( {- {\int_0}^{L_1^{+1} \dot{X}_1 \, L_2^{-1} \dot{Y}_2} {\, 
\log(1+t) \, \over \, t \,} \; dt}\,{\bigg|}_{\hbar=0} \right) 
\right) =   \hfill  \cr   
   = \exp \! \left( {\hbox{\rm ad}}_{\{\ ,\ \}} \left(
- {\int_0}^{z_1^{-1} x_1 \, z_2^{+1} y_2} {\, \log(1+t) \, \over 
\, t \,} \; dt \right) \right) =  \cr   
   \hfill   = \exp \! \left( \mu \left( - {\, \log \left( 1 +
z_1^{-1} x_1 \, z_2^{+1} y_2 \right) \, \over \, z_1^{-1} x_1 \, 
z_2^{+1} y_2 \,} \right) \circ {\hbox{\rm ad}}_{\{\ ,\ \}} \left( 
z_1^{-1} x_1 \, z_2^{+1} y_2 \right) \right)  \cr }  $$   where  $ 
\, \mu(S) \, $  denotes the operator of left multiplication by  $ 
\, S \in F \big[\big[ {\gersl_2}^{\!*} \oplus {\gersl_2}^{\!*} 
\big]\big] \, $.  Indeed   
  $$  \displaylines{
   \hbox{\rm ad}_{\{ \ , \ \}} \left( - {\int_0}^{z_1^{-1} x_1 \,
z_2^{+1} y_2} {\, \log(1+t) \, \over \, t \,} \; dt \right) (\ell) 
= \hbox{\rm ad}_{\{ \ , \ \}} \left( - \sum_{n>0} \big( z_1^{-1} x_1 
\, z_2^{+1} y_2 \big) \Big/ n^2 \right) (\ell) =   \hfill  \cr   
   = - \sum_{n>0} {1 \over n^2} \cdot \left\{ {\big( z_1^{-1} x_1
\, z_2^{+1} y_2 \big)}^n, \ell \right\} = - \sum_{n>0} {\,1\, 
\over \,n^2\,} \, n {\big( z_1^{-1} x_1 \, z_2^{+1} y_2 
\big)}^{n-1} \cdot \big\{ z_1^{-1} x_1 \, z_2^{+1} y_2 \, , \ell 
\, \big\} =  \cr   
   \hfill   = - \sum_{n>0} {{\big( z_1^{-1} x_1 \, z_2^{+1} y_2
\big)}^{n-1} \over n} \cdot \big\{ z_1^{-1} x_1 \, z_2^{+1} y_2 \, 
, \ell \big\} = - {\log \left( 1 + z_1^{-1} x_1 \, z_2^{+1} y_2 
\right) \over \, z_1^{-1} x_1 \, z_2^{+1} y_2 \,} \cdot \big\{ 
z_1^{-1} x_1 \, z_2^{+1} y_2 \, , \ell \big\}  \cr }  $$   
(because of Jacobi identity:  $ \, \{ \, \cdot \, , \ell\} = - 
\hbox{\rm ad}_{\{ \ , \ \}}(\ell)  \, $  is a derivation!).  Second, 
again by the Jacobi identity and the commutation relation  $ \, 
z_1^{-1} x_1 \cdot z_2^{+1} y_2 = z_2^{+1} y_2 \cdot z_1^{-1} x_1 
\, $  we get  
  $$  \displaylines{ 
   \mu \left( - {\, \log \left( 1 + z_1^{-1} x_1 \, z_2^{+1} y_2
\right) \, \over \, z_1^{-1} x_1 \, z_2^{+1} y_2 \,} \right) \circ 
{\hbox{\rm ad}}_{\{\ ,\ \}} \left( z_1^{-1} x_1 \, z_2^{+1} y_2 
\right) \; =   \hfill  \cr   
   = \; \mu \left( - {\, \log \left( 1 + z_1^{-1} x_1 \, z_2^{+1}
y_2 \right) \, \over \, z_1^{-1} x_1 \, z_2^{+1} y_2 \,} \right)
\circ {\hbox{\rm ad}}_{\{\ ,\ \}} \left( z_1^{-1} x_1 \right)
\circ \mu \left( z_2^{+1} y_2 \right) \; +  \cr   
   \hfill   + \; \mu \left( - {\, \log \left( 1 + z_1^{-1} x_1 \,
z_2^{+1} y_2 \right) \, \over \, z_1^{-1} x_1 \, z_2^{+1} y_2 \,} 
\right) \circ \mu \left( z_1^{-1} x_1 \right) \circ {\hbox{\rm 
ad}}_{\{\ ,\ \}} \left( z_2^{+1} y_2 \right) \; =  \cr   
   = \; \mu \left( - {\, \log \left( 1 + z_1^{-1} x_1 \, z_2^{+1}
y_2 \right) \, \over \, z_1^{-1} x_1 \,} \right) \circ {\hbox{\rm 
ad}}_{\{\ ,\ \}} \left( z_1^{-1} x_1 \right) \; +  \hfill  \cr   
   \hfill   + \; \mu \left( - {\, \log \left( 1 + z_2^{+1} y_2
\right) \, \over \, z_1^{-1} x_1 \, z_2^{+1} y_2 \,} \right) \circ 
{\hbox{\rm ad}}_{\{\ ,\ \}} \left( z_2^{+1} y_2 \right) \; ; 
\cr }  $$   
the two summands above are mutually commuting operators   
--- thanks to the commutation relation  $ \, z_1^{-1} x_1 \cdot 
z_2^{+1} y_2 = z_2^{+1} y_2 \cdot z_1^{-1} x_1 $  ---   so when we 
take the exponential we get   
  $$  \displaylines{ 
   \R_{{}_{G\!H}}^{(1)} = \exp \! \left( \mu \left( - {\, \log \left(
1 + z_1^{-1} x_1 \, z_2^{+1} y_2 \right) \, \over \, z_1^{-1} x_1 
\, z_2^{+1} y_2 \,} \right) \circ {\hbox{\rm ad}}_{\{\ ,\ \}} 
\left( z_1^{-1} x_1 \, z_2^{+1} y_2 \right) \right) =   \hfill  
\cr   
   = \exp \! \Bigg( \mu \left( - {\, \log \left( 1 + z_1^{-1} x_1
\, z_2^{+1} y_2 \right) \, \over \, z_1^{-1} x_1 \,} \right) \circ 
{\hbox{\rm ad}}_{\{\ ,\ \}} \left( z_1^{-1} x_1 \right) \; +  \cr   
   \hfill   + \; \mu \left( - {\, \log \left( 1 + z_1^{-1} x_1 \,
z_2^{+1} y_2 \right) \, \over \, z_2^{+1} y_2 \,} \right) \circ 
{\hbox{\rm ad}}_{\{\ ,\ \}} \left( z_2^{+1} y_2 \right) \Bigg) =  
\cr 
   = \exp \! \left( \mu \left( - {\, \log \left( 1 + z_1^{-1} x_1
\, z_2^{+1} y_2 \right) \, \over \, z_1^{-1} x_1 \,} \right) \circ 
{\hbox{\rm ad}}_{\{\ ,\ \}} \left( z_1^{-1} x_1 \right) \right) \; 
\circ   \hfill  \cr   
   \hfill  \circ \; \exp \! \left( \mu \left( - {\, \log \left( 1
+ z_1^{-1} x_1 \, z_2^{+1} y_2 \right) \, \over \, z_2^{+1} y_2 
\,} \right) \circ {\hbox{\rm ad}}_{\{\ ,\ \}} \left( z_2^{+1} y_2 
\right) \right) \, ;  \cr }  $$   in a nutshell, we have   
  $$  \R_{{}_{G\!H}}^{(1)} = \exp \! \left( \Cal{E}_1 \right)
\circ \exp \! \left( \Cal{F}_2 \right)  \hskip7pt  \hbox{with} 
\hskip6pt  \Bigg\{ \hbox{ $ \eqalign{ 
   \Cal{E}_1  &  := \mu \left( - \log \left( \nabla^2 \right)
\Big/ z_1^{\!-1} x_1 \right) \circ {\hbox{\rm ad}}_{\{\ ,\ \}} 
\left( z_1^{\!-1} x_1 \right)  \cr   
   \Cal{F}_2  &  := \mu \left( - \log \left(  \nabla^2 \right)
\Big/ z_2^{\!+1} y_2 \right) \circ {\hbox{\rm ad}}_{\{\ ,\ \}} 
\left( z_2^{\!+1} y_2 \right)  \cr } $ }   \hskip4,7pt (4.16)  $$   
where  $ \; \nabla := {\big( 1 + z_1^{\!-1} x_1 \, z_2^{\!+1} y_2 
\big)}^{1/2} \, $.  We proceed now with computations.   
                                            \par  
   Since  $ \, \{s_1,r_2\} = 0 \, $  for all  $ \, s, r
\in F\big[\big[{\gersl_2}^{\!*}\big]\big] $,  we have  $ \, 
\Cal{E}_1(r_2) = 0 \, $  for all  $ \, r \in F\big[\big[ 
{\gersl_2}^{\!*}\big]\big] $,  \, so   
  $$  \exp \! \left( \Cal{E}_1 \right)(x_2) = x_2 \; ,  \qquad 
\exp \! \left( \Cal{E}_1 \right)\big(z_2^{\!\pm 1}\big) = 
z_2^{\!\pm 1} \; ,  \qquad  \exp \! \left( \Cal{E}_1 \right)(y_2) 
= y_2 \; .   \eqno (4.17)  $$   
   \indent   Now for the rest!  We have to compute things like  $ \,
\{s_1,r_1\} \, $,  \, so for simplicity we shall drop the 
subscript 1 throughout.   
                                            \par  
   For the operator  $ \; \hbox{\rm ad}_{\{\ ,\ \}}\big(z^{-1}x\big)
\; $  (a derivation!) we have the formul\ae{}   
  $$  \displaylines{  
   \qquad  \big\{ z^{-1} x, x \big\} = \big\{ z^{-1}, x \big\} \cdot x
+ z^{-1} \cdot \big\{ x, x \big\} = \big\{ z^{-1}, x \big\} \cdot 
x = \big( z^{-1} x \big/ 2 \big) \cdot x   \hfill  \cr   
   \qquad  \big\{ z^{-1} x, z^{\pm 1} \big\} = \big\{ z^{-1}, z^{\pm 1}
\big\} \cdot x + z^{-1} \cdot \big\{ x, z^{\pm 1} \big\} = z^{-1} 
\cdot \big\{ x, z^{\pm 1} \big\} = \pm \big( z^{-1} x \big/ 2 
\big) \cdot z^{\pm 1}   \hfill   \cr   
   \qquad  \big\{ z^{-1} x, y \big\} = \big\{ z^{-1}, y \big\} \cdot x
+ z^{-1} \cdot \big\{ x, y \big\} = - \big( z^{-1} \big/ 2 \big) 
\cdot y \cdot x + z^{-1} \cdot \big( z^{-2} - z^{+2} \big) =   
\hfill  \cr   
   \hfill   = - \big( z^{-1} x \big/ 2 \big) \cdot y + z^{-3} - z^{+1} 
\cr   
   \qquad  \big\{ z^{-1} x, z^{+1} y \big\} = \big\{ z^{-1} x, z^{+1}
\big\} \cdot y + z^{+1} \cdot \big\{ z^{-1} x, y \big\} = z^{-2} - 
z^{+2}   \hfill  \cr }  $$   
   \indent   Then for  $ \; \exp \! \left( \Cal{E}_1 \right) = \exp \!
\left( \mu \left( - \log \left( \nabla^2 \right) \Big/ z_1^{\!-1} 
x_1 \right) \circ {\hbox{\rm ad}}_{\{\ ,\ \}} \left( z_1^{\!-1} 
x_1 \right) \right) \; $  we have   
  $$  \hbox{ $ \eqalign{ 
   \exp \! \left( \Cal{E}_1 \right)  &  (x_1) = \exp \Big( \! - \log
\! \big( \nabla^2 \big) \Big/ 2 \Big) \cdot x_1  = \exp \big(\!-\! 
\log(\nabla) \big) \cdot x_1 = x_1 \cdot \nabla^{-1}  \cr   
   \exp \! \left( \Cal{E}_1 \right)  &  \big(z_1^{\!\pm 1}\big)
= \exp \Big( \! \mp \log \! \big( \nabla^2 \big) \Big/ 2 \Big) 
\cdot z_1^{\!\pm 1} = \exp \big( \! \mp \! \log(\nabla) \big) 
\cdot z_1^{\!\pm 1} = z_1^{\!\pm 1} \cdot \nabla^{\mp 1}  \cr   
   \exp \! \left( \Cal{E}_1 \right)  &  (y_1) = y_1 \cdot
\nabla^{+1} + y_2 \cdot z_2^{\!+1} z_1^{\!-3} \cdot \nabla^{+1} - 
y_2 \cdot z_2^{\!+1} z_1^{\!+1} \cdot \nabla^{-1}  \cr } $ }  
\eqno (4.18)  $$   
where the latter identity is computed (since  $ \, \exp \! \left(
\Cal{E}_1 \right) \, $  is an automorphism!) as follows:    
  $$  \displaylines{ 
   \exp \! \left( \Cal{E}_1 \right)(y_1) = \exp \! \left( \Cal{E}_1
\right)\big( z_1^{\!-1} \cdot z_1^{\!+1} y_1 \big) = \exp \! 
\left( \Cal{E}_1 \right)\big( z_1^{\!-1} \big) \cdot \exp \! 
\left( \Cal{E}_1 \right) \big( z_1^{\!+1} y_1 \big) =   \hfill  
\cr  
   = z_1^{\!-1} \cdot \nabla^{+1} \cdot \bigg( \, z_1^{\!+1} y_1
+ \sum_{n>1} {\,1\, \over \,n!\,} \cdot {\Big( \! - \! \log \left( 
\nabla^2 \right) \! \Big/ z_1^{\!-1} x_1 \Big)}^n \cdot 
{\Cal{E}_1}^{\! n-1} \Big( z_1^{\!-2} - z_1^{\!+2} \Big) \bigg) =  
\cr   
   = z_1^{\!-1} \cdot \nabla^{+1} \cdot \bigg( \, z_1^{\!+1} y_1 +
\sum_{n>0} {\,1\, \over \,n!\,} \cdot {\Big( \! - \! \log \left( 
\nabla^2 \right) \! \Big/ z_1^{\!-1} x_1 \Big)}^n \cdot {\big( - 
z_1^{\!-1} x_1 \big)}^{n-1} \cdot z_1^{\!-2} \; -   \hfill  \cr   
   \hfill   - \; \sum_{n>0} {\,1\, \over \,n!\,} \cdot {\Big( \!
- \! \log \left( \nabla^2 \right) \! \Big/ z_1^{\!-1} x_1 \Big)}^n 
\cdot {\big( + z_1^{\!-1} x_1 \big)}^{n-1} \cdot z_1^{\!+2} \Big) 
\bigg) =  \cr   
   = z_1^{\!-1} \cdot \nabla^{+1} \cdot \bigg( \, z_1^{\!+1} y_1 + 
{\; \exp \big( \! + \log \big( \nabla^2 \big) \big) - 1 \; \over 
\; + z_1^{\!-1} x_1 \;} \cdot z_1^{\!-2} - {\; \exp \big( \! - 
\log \big( \nabla^2 \big) \big) - 1 \; \over \; - z_1^{\!-1} x_1 
\;} \cdot z_1^{\!+2} \bigg) =  \cr   
   \hfill   = z_1^{\!-1} \cdot \nabla^{+1} \cdot \bigg( \, z_1^{\!+1}
y_1 + {\; \nabla^{+2} - 1 \; \over \; + z_1^{\!-1} x_1 \;} \cdot 
z_1^{\!-2} - {\; \nabla^{-2} - 1 \; \over \; - z_1^{\!-1} x_1 \;} 
\cdot z_1^{\!+2} \bigg) =  \cr   
   = z_1^{\!-1} \cdot \nabla^{+1} \cdot \bigg( \, z_1^{\!+1} y_1 + 
y_2 \cdot z_2^{\!+1} z_1^{\!-2} - y_2 \cdot z_2^{\!+1} z_1^{\!+2} 
\cdot \nabla^{-2} \bigg) =   \hfill  \cr   
   \hfill   = y_1 \cdot \nabla^{+1} + y_2 \cdot z_2^{\!+1} z_1^{\!-3}
\cdot \nabla^{+1} - y_2 \cdot z_2^{\!+1} z_1^{\!+1} \cdot 
\nabla^{-1} \; .  \cr }  $$   
   \indent   Now for  $ \; \exp \! \left( \Cal{F}_2 \right) \, $. 
Again, since  $ \, \{s_1,r_2\} = 0 \, $  for all  $ \, s, r \in 
F\big[\big[{\gersl_2}^{\!*}\big]\big] \, $  we have  $ \, 
\Cal{F}_2 (s_1) = 0 \, $  for all  $ \, s \in F \big[\big[ 
{\gersl_2}^{\!*} \big]\big] $,  \; so   
  $$  \exp \! \left( \Cal{F}_2 \right)(x_1) = x_1 \; ,  \qquad 
\exp \! \left( \Cal{F}_2 \right)\big(z_1^{\!\pm 1}\big) = 
z_1^{\!\pm 1} \; ,  \qquad  \exp \! \left( \Cal{F}_2 \right)(y_1) 
= y_1 \; .   \eqno (4.19)  $$   
   As for the rest, we can base upon the previous results, as
follows.  First, we note that there is a  {\sl Poisson algebra\/} 
automorphism  
  $$  \Phi \, \colon \, F\big[\big[{\gersl_2}^{\!*}\big]\big]
\; {\buildrel \cong \over {\llongrightarrow}} \; 
F\big[\big[{\gersl_2}^{\!*}\big]\big] \; ,  \qquad  x \mapsto y \, 
,  \quad z^{\pm 1} \mapsto z^{\mp 1} \, ,  \quad  y \mapsto x  $$  
such that  $ \, \Phi^{-1} = \Phi \, $  (and which also restrict to 
$ \, F \big[ {}_s{{SL}_2}^{\!*} \big] \, $  and to  $ \, F \big[ 
{}_a{{SL}_2}^{\!*}\big] \, $).  Then we have immediately from 
definitions that  $ \; (\Phi \otimes \Phi) \big( \Cal{E}_1 \big) 
(s \otimes r) = \sigma \Big( \Cal{F}_2 \big( \Phi(r) \otimes 
\Phi(s) \Big) \; $  for all  $ \; s, r \in F \big[\big[ 
{\gersl_2}^{\!*} \big]\big] \; $  (with  $ \sigma $  as in \S 
1.8), whence in particular we argue   
  $$  \Cal{F}_2(s_2) = \sigma \Big( \Phi^{\otimes 2} \Big(
\Cal{E}_1 \big( \Phi^{-1}(s_1) \big) \! \Big) \Big) = \sigma \Big( 
\Phi^{\otimes 2} \Big( \Cal{E}_1 \big( \Phi(s_1) \big) \! \Big) 
\Big)  \qquad  \forall \;\; s \in F \big[\big[ {\gersl_2}^{\!*} 
\big]\big]  $$   
   and so
  $$  \exp\big(\Cal{F}_2\big)(s_2) = \sigma \Big( \Phi^{\otimes 2}
\Big( \exp\big(\Cal{E}_1\big) \big( \Phi(s_1) \big) \! \Big) \Big) 
\qquad  \forall \;\; s \in F \big[\big[ {\gersl_2}^{\!*} 
\big]\big] \; .  $$   
Using this and formul\ae{}  (4.18) we eventually get   
  $$  \hbox{ $ \eqalign{ 
   \exp \! \left( \Cal{F}_2 \right)  &  (x_2) = x_2 \cdot \nabla^{+1}
+ x_1 \cdot z_1^{\!-1} z_2^{\!+3} \cdot \nabla^{+1} - x_1 \cdot 
z_1^{\!-1} z_2^{\!-1} \cdot \nabla^{-1}  \cr   
   \exp \! \left( \Cal{F}_2 \right)  &  \big(z_2^{\!\pm 1}\big) =
z_2^{\!\pm 1} \cdot \nabla^{\pm 1}  \cr   
   \exp \! \left( \Cal{F}_2 \right)  &  (y_1) =
y_2 \cdot \nabla^{-1} \; .  \cr } $ }   \eqno (4.20)  $$   
   \indent   Formul\ae{}  (4.16--20) give us a complete description
of  $ \, \R^{(1)}_{{}_{G\!H}} \, $:  \, to summarize, it is given 
by   
  $$  \hbox{ $ \eqalign{ 
   \R^{(1)}_{{}_{G\!H}} \big( x_1 \big) = 
x_1 \cdot \nabla^{-1} \; ,  \qquad \qquad  \R^{(1)}_{{}_{G\!H}} 
\big( {z_1}^{\!\pm 1} \big) = {z_1}^{\!\pm 1} \cdot \nabla^{\mp 1} 
\qquad  \cr  
   \R^{(1)}_{{}_{G\!H}} \big( y_1 \big) = 
y_1 \cdot \nabla^{+1}  + y_2 \cdot {z_2}^{\!+1} {z_1}^{\!-3} \cdot 
\nabla^{+1} - y_2 \cdot {z_2}^{\!+1} {z_1}^{\!+1} \cdot 
\nabla^{-1} \cr   
   \R^{(1)}_{{}_{G\!H}} \big( x_2 \big) = 
x_2 \cdot \nabla^{+1}  + x_1 \cdot {z_1}^{\!-1} {z_2}^{\!+3} \cdot 
\nabla^{+1} - x_1 \cdot {z_1}^{\!-1} {z_2}^{\!-1} \cdot 
\nabla^{-1} \cr   
   \R^{(1)}_{{}_{G\!H}} \big( {z_2}^{\!\pm 1} \big) = 
{z_2}^{\!\pm 1} \cdot \nabla^{\pm 1} \; ,  \qquad \qquad 
\R^{(1)}_{{}_{G\!H}} \big( y_2 \big) =  y_2 \cdot \Delta^{-1} 
\qquad  \cr } $ }   \eqno (4.21)  $$   
   \indent   Finally, composing with  $ \, \R^{(0)}_{{}_{G\!H}} \, $  
--- see (4.15) ---   we find at last   
  $$  \hbox{ $ \eqalign{ 
   \R_{{}_{G\!H}} \big( x_1 \big) = 
x_1 \cdot {z_2}^{\! -2} \cdot \Theta^{-1} \; ,  \qquad \qquad 
\R_{{}_{G\!H}} \big( {z_1}^{\!\pm 1} \big) =  {z_1}^{\!\pm 1} 
\cdot \Theta^{\mp 1}  \qquad  \cr  
   \R_{{}_{G\!H}} \big( y_1 \big) = 
y_1 \cdot {z_2}^{\!+2} \cdot \Theta^{+1}  + y_2 \cdot {z_2}^{\!+1} 
{z_1}^{\!-1} \cdot \Theta^{+1} - y_2 \cdot {z_2}^{\!+1} 
{z_1}^{\!+3} \cdot \Theta^{-1}  \cr   
   \R_{{}_{G\!H}} \big( x_2 \big) = 
x_2 \cdot {z_1}^{\!-2} \cdot \Theta^{+1}  + x_1 \cdot {z_1}^{\!-1} 
{z_2}^{\!+1} \cdot \Theta^{+1} - x_1 \cdot {z_1}^{\!-1} 
{z_2}^{\!-3} \cdot \Theta^{-1}  \cr   
   \R_{{}_{G\!H}} \big( {z_2}^{\!\pm 1} \big) = 
{z_2}^{\!\pm 1} \cdot \Theta^{\pm 1} \; ,  \qquad \qquad 
\R_{{}_{G\!H}} \big( y_2 \big) =  y_2 \cdot {z_1}^{\! +2} \cdot 
\Theta^{-1}  \qquad  \cr } $ }   \eqno (4.22)  $$   for  $ \; 
\R_{{}_{G\!H}} = \R^{(0)}_{{}_{G\!H}} \circ \R^{(1)}_{{}_{G\!H}} 
\; $  (see (4.14)), with  $ \; \Theta := {\big( 1 + x_1 z_1^{\!+1} 
z_2^{\!-1} y_2 \big)}^{1/2} = \R^{(0)}_{{}_{G\!H}} \big( \nabla 
\big) \, $.   
                                                \par
   Therefore, just comparing (4.22) with (4.8) we get as an outcome
the main result of this section:  

\vskip9pt   

\proclaim{Theorem 4.6}  The braidings  $ \; \R_{{}_{W\!X}} \, $ 
and  $ \; \R_{{}_{G\!H}} \, $  for  $ \, \gerg = \gersl_2({\Bbb C})
\, $  do coincide.  In other words, the answer to the "Question"
in \S 4.1 is positive for  $ \, \gerg = \gersl_2({\Bbb C}) \, $.  
\qed   
\endproclaim

%%  
%%  \vskip9pt   
%%  
%%  
%%  {\bf Remarques finales:} \, dans la  \S 1.8{\it (b)}  j'ecris
%%  {\it "The algebra is called  {\sl rigid},  and the braiding
%%  operator  {\sl unitary},  if in addition [...]"}.  Bon, cela
%%  est une "tentative terminology" que je viens d'inventer: demande
%%  s'il-te-plait \`a Turaev s'il est au courant d'une quelque
%%  terminologie dej\`a etablie pour cette notion (\`a savoir,
%%  la sous-classe des toutes les alg\`ebres de Hopf tress\'ees
%%  qui correspndent aux alg\`ebres de Hopf triangulaires).
%%  
%%  

\vskip15pt

\enddocument

\bye